\documentclass[a4paper,10pt]{article}
\usepackage{amsmath, amssymb, algorithm, algorithmic}
\usepackage{graphicx}
\usepackage{psfrag}

\newcommand{\dt}[1]{\dot{\theta_#1}}

\newcommand{\R}{\mathbb{R}}
\newcommand{\C}{\mathbb{C}}
\newcommand{\Z}{\mathbb{Z}}
\newcommand{\NewtonPolygon}[1]{\Pi(#1)}
\newcommand{\supp}[1]{\mathrm{Support}(#1)}
\newcommand{\coef}[1]{\mathrm{Coef}(#1)}
\newcommand{\NP}[1]{\Pi(#1)}

\newcommand{\GName}{G}
\newcommand{\FName}{F}
\newcommand{\FNameCoeff}{{\bf F}}
\newcommand{\HName}{H}
\newcommand{\KName}{K}

\newcommand{\p}{{\bf p}}
\newcommand{\z}{{\bf z}}

\newtheorem{Definition}{Definition}

\newtheorem{Proof}{Proof}
\newtheorem{Theorem}{Theorem}

\newtheorem{Example}{Example}
\newtheorem{Remark}{Remark}

\title{Dynamic balancing of planar mechanisms using toric geometry\footnote{B. Moore and J. Schicho were partially supported by the Austrian Science
Fund (FWF) under the SFB grant F1303. B. Moore was also supported by the Fonds Qu{\'e}b{\'e}cois de la Recherche sur la Nature et les Technologies (FQRNT). C. Gosselin was partially supported by the Natural Sciences and Engineering Research Council of Canada (NSERC).}}
\author{
\small Brian Moore, Josef Schicho\\
\small Johann Radon Institute for Computational and Applied Mathematics\\
\small Austrian Academy of Science, A-4040, Linz, Austria \\
\scriptsize \{brian.moore, josef.schicho\}@ricam.oeaw.ac.at\\ 
\\
\large
\small Cl\'ement M. Gosselin\\
\small D\'epartement de G\'enie M\'ecanique, Universit\'e Laval, Qu\'ebec\\
\small Qu\'ebec, G1K 7P4, Canada, \\
\scriptsize gosselin@gmc.ulaval.ca\\
\date{}
}
\normalsize

\begin{document}
\maketitle

\begin{abstract}
In this paper, a new method to determine the complete set of dynamically 
balanced planar four-bar mechanims is presented. 
Using complex variables to model the kinematics of the mechanism, 
the dynamic balancing constraints are written as
algebraic equations over complex variables and joint angular velocities.
After elimination of the joint angular velocity variables, 
the problem is formulated as a problem of factorization of Laurent polynomials.
Using toric polynomial division, necessary and sufficient conditions for 
dynamic balancing of planar four-bar mechanisms are derived.
\end{abstract}



 \section{Introduction}
%
%
%


Statically and dynamically balanced mechanisms are highly desirable for many engineering applications since
they do not apply any forces or moments at their base, for arbitrary motion trajectories.
This concept is used in mechanism design in order to reduce fatigue, vibrations and wear.
Static and dynamic balancing can also be used in more advanced applications such as, for instance, 
the design of more efficient flight simulators \cite{UphoffGosselinLaliberte:2000}, or in the
design of compensation mechanisms for telescopes.  
Additionally, dynamic balancing is very attractive for space applications since the
reaction forces induced at the base of space manipulators or mechanisms are one of the reasons
why the latter are constrained to move very slowly \cite{Yoshida:2001}.

Several approaches can be used to balance mechanisms (see for instance \cite{Bagci:1982, ArakelianSmith:1999}).
In general, complete balancing requires the integration of additional mechanical components in
the design of a mechanism, such as counterweights and counterrotations \cite{Berkof_Lowen:1971}. However, for
some simple architectures, it is sometimes possible to design dynamically balanced mechanisms
by an appropriate choice of the design parameters without introducing additional
linkages or counterrotations \cite{GosselinVollmerCoteWu:2004}.  Even though these mechanisms
do not include counterrotations, they satisfy all conditions for balancing, namely, their
centre of mass remains stationary (static balancing) and
their total angular momentum vanishes (dynamic balancing), for arbitrary trajectories.

Although families of dynamically balanced four-bar mechanisms were presented in
\cite{GosselinVollmerCoteWu:2004}, no proof was given on the existence of other
possible solutions. 
In this paper, the aim is to derive all possible sets of design parameters for which a planar four-bar mechanism
is dynamically balanced. 
The problem was first addressed by Berkof and Lowen \cite{Berkof_Lowen:1969}, 
who provided conditions for static balancing in terms of the 
design parameters when the geometric parameters are sufficiently generic.
A non-generic solution was then found in \cite{Gosselin:1997}.
In \cite{MooreSchichoGosselin:2007}, a complete list of
statically balanced planar four-bar mechanisms was given. The problem of dynamic
balancing was also addressed in \cite{RicardGosselin:2000}, where special cases
that do not require external counterrotations were first revealed. In the latter reference,
it can be observed that the dynamic balancing problem leads to a rather complicated system of 
algebraic equations.

A generic approach to investigate such parametric polynomial systems has been
recently proposed by Lazard and Rouillier \cite{LazardRouillier:2007}.
In \cite{MooreSchichoGosselin:2007}, it was shown that the
problem can be simplified if one uses complex variables to model the
angles in the configuration space, together with some results in
algebraic geometry (Ostrowski's theorem ~\cite{Ostrowski:1921}).
However, the application of the above techniques to the problem of
dynamic balancing leads to a cumbersome case by case analysis. In order to
simplify this analysis, a concept of division for
Laurent polynomials is introduced here. The method is
related to \cite{SalemGaoLauder:2004}, who use similar ideas for
factoring polynomials by taking advantage of the special shape
of their Newton polyhedra.



This paper is organized as follows. First, we describe how to derive the kinematic, static and dynamic equations
using complex number representations. This leads to a set of algebraic equations in two complex variables where
the coefficients are expressions in terms of the design parameters.
This section can be skipped if one is interested only in the mathematical aspects.
Then, we introduce the concept of division for Laurent polynomials and
provide an algorithm for computing this division.
Then, we use Laurent polynomial division to eliminate the dependent variables
describing the configuration space. Finally, we solve the system
in the remaining design parameters and give an overview of the
various subcases.


\section{Model}


\subsection{Representation of planar four-bar mechanisms}

\begin{figure}
	\centering
        \psfrag{l1}{$l_1$}
        \psfrag{l2}{$l_3$}
        \psfrag{l3}{$l_2$}
        \psfrag{d}{$d$}
        \psfrag{m1}{$m_1$}
        \psfrag{m2}{$m_3$}
        \psfrag{m3}{$m_2$}
        \psfrag{t1}{$\theta_1(t)$}
        \psfrag{t2}{$\theta_3(t)$}
        \psfrag{t3}{$\theta_2(t)$}
        \psfrag{k1}{$\psi_1$}
        \psfrag{k2}{$\psi_3$}
        \psfrag{k3}{$\psi_2$}
        \psfrag{r1}{$r_1$}
        \psfrag{r2}{$r_3$}
        \psfrag{r3}{$r_2$}
        \psfrag{X}{}
        \psfrag{Y}{}
	\includegraphics[height=5.5cm]{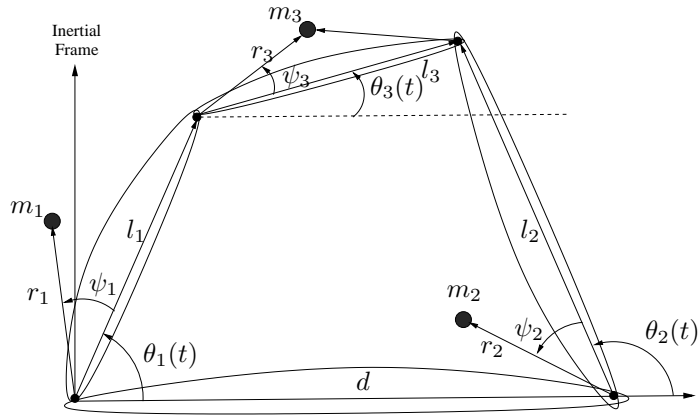}
	\caption{Four-bar mechanism.}
	\label{fig:bar4}
\end{figure}

A planar four-bar mechanism is shown in Figure \ref{fig:bar4}.
It consists of four links: the base of length $d$ which is fixed, and three moveable links of length $l_1, l_2, l_3$ respectively.
We assume that all link lengths are strictly positive.
Since the base is fixed, the mass properties of the base has no influence on the equations and will therefore be ignored.
Each of the three moveable links has a mass $m_i$, a centre of mass whose position is defined by $r_i$ and $\psi_i$ and a moment of inertia $I_i$.
The design of planar four-bar mechanisms consists in choosing the 16 design parameters (Table \ref{tab:designParam}).

\begin{table}
\begin{tabular}{|c|c||c|}
\hline
\hspace*{14mm}Type\hspace*{14mm} & \hspace*{14mm}\hspace*{14mm} & \hspace*{14mm}Parameters\hspace*{14mm}\\
\hline
\hline
&&\\
{\bf Kinematic} & Length & $l_1, l_2, l_3, d$ \\
&&\\
\hline
&&\\
{\bf Static} & Mass & $m_1, m_2, m_3$ \\
&&\\
              & Centre of mass & $r_1, \psi_1, r_2, \psi_2, r_3, \psi_3$ \\
&&\\
\hline
&&\\
{\bf Dynamic} & Inertia & $I_1, I_2, I_3$ \\
&&\\
\hline
\end{tabular}
\caption{Design parameters for the planar four-bar mechanisms.}
\label{tab:designParam}
\end{table}

\begin{figure}
	\begin{tabular}{cc}
 		\includegraphics[width=5.5cm]{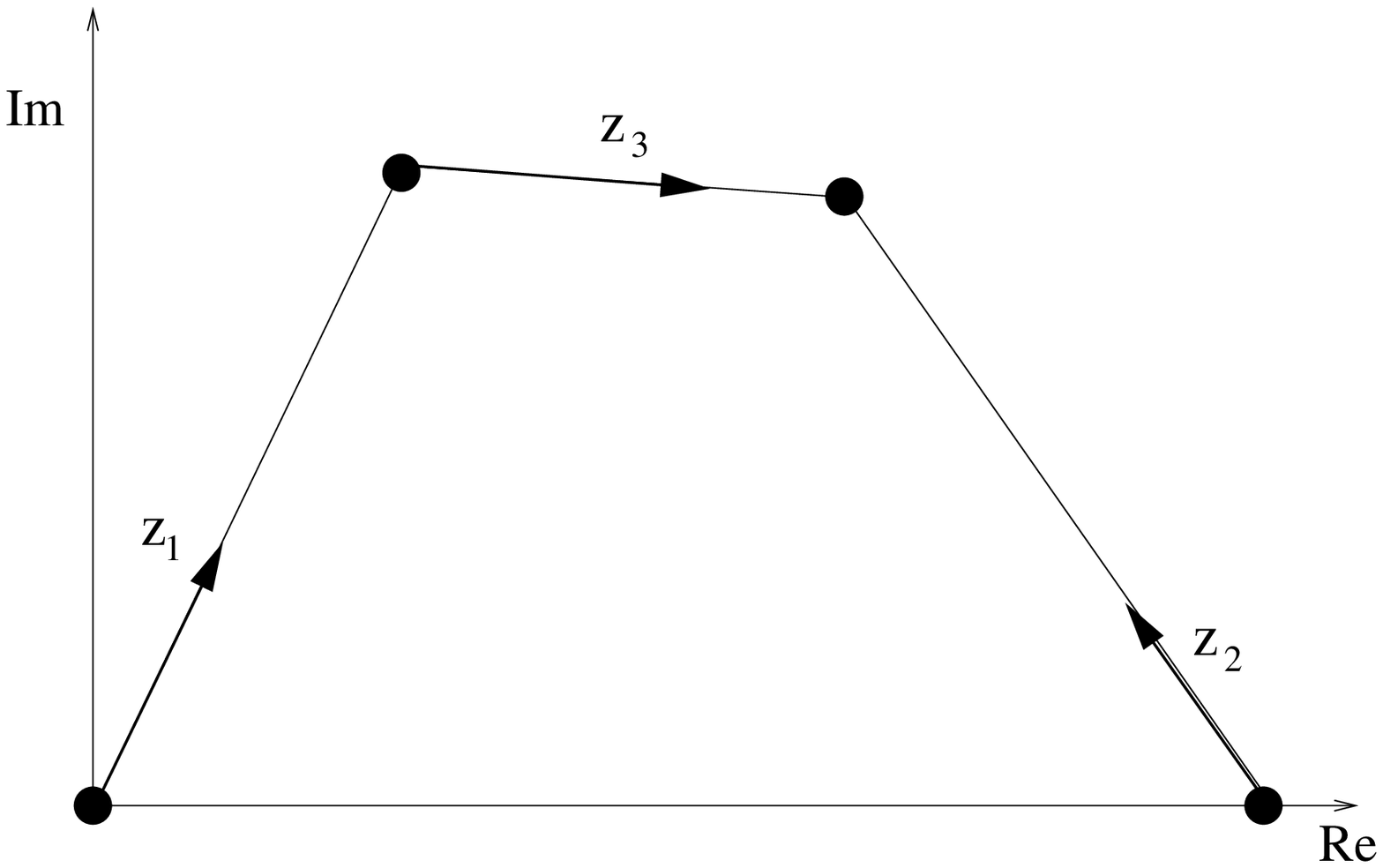} &
		\includegraphics[width=5.5cm]{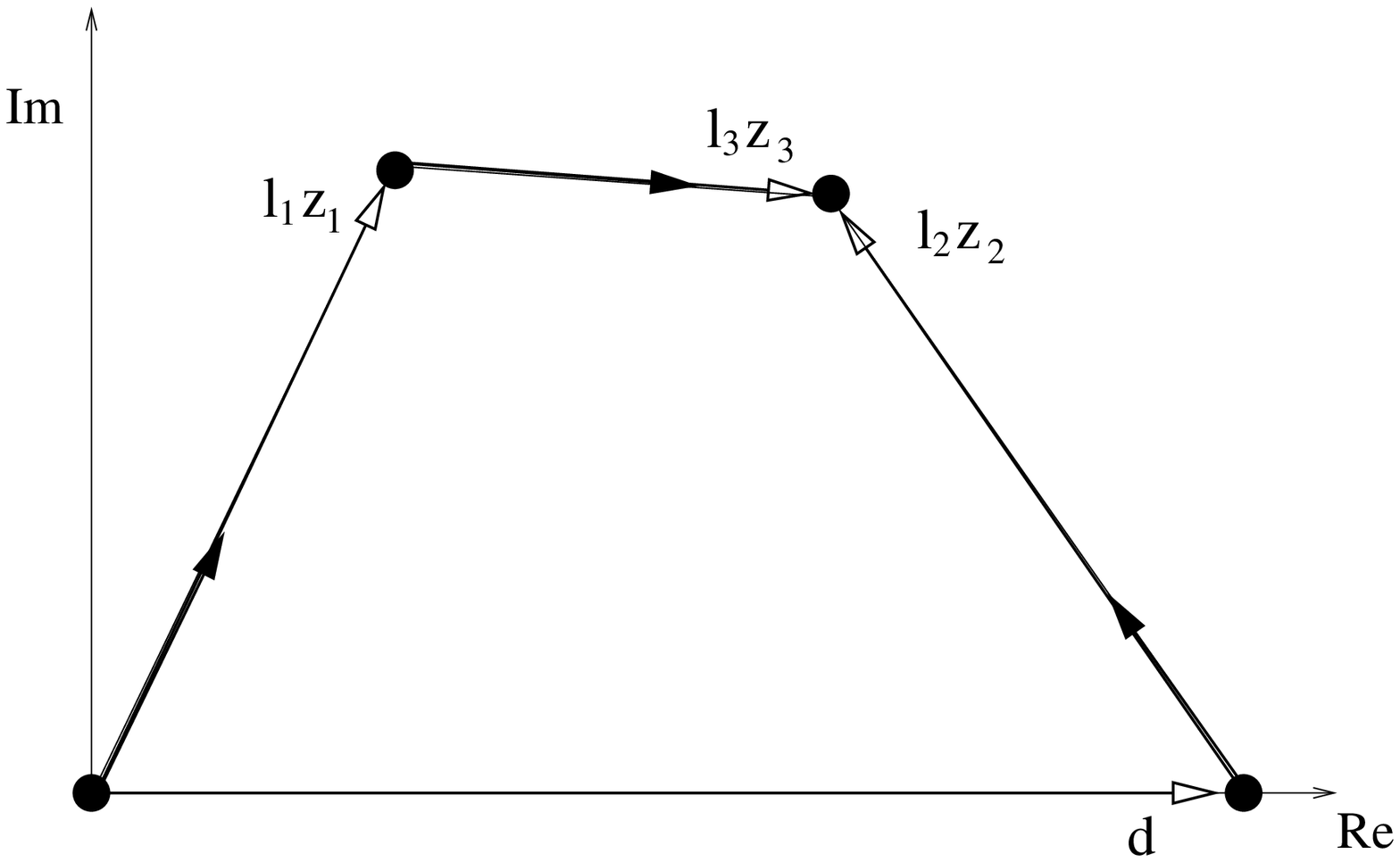} 
	\end{tabular} 
\caption{Complex representation for the kinematics.}
\label{fig:bar4Representation}
\end{figure}

\begin{figure}
\begin{tabular}{cc}
 	\includegraphics[width=5.5cm]{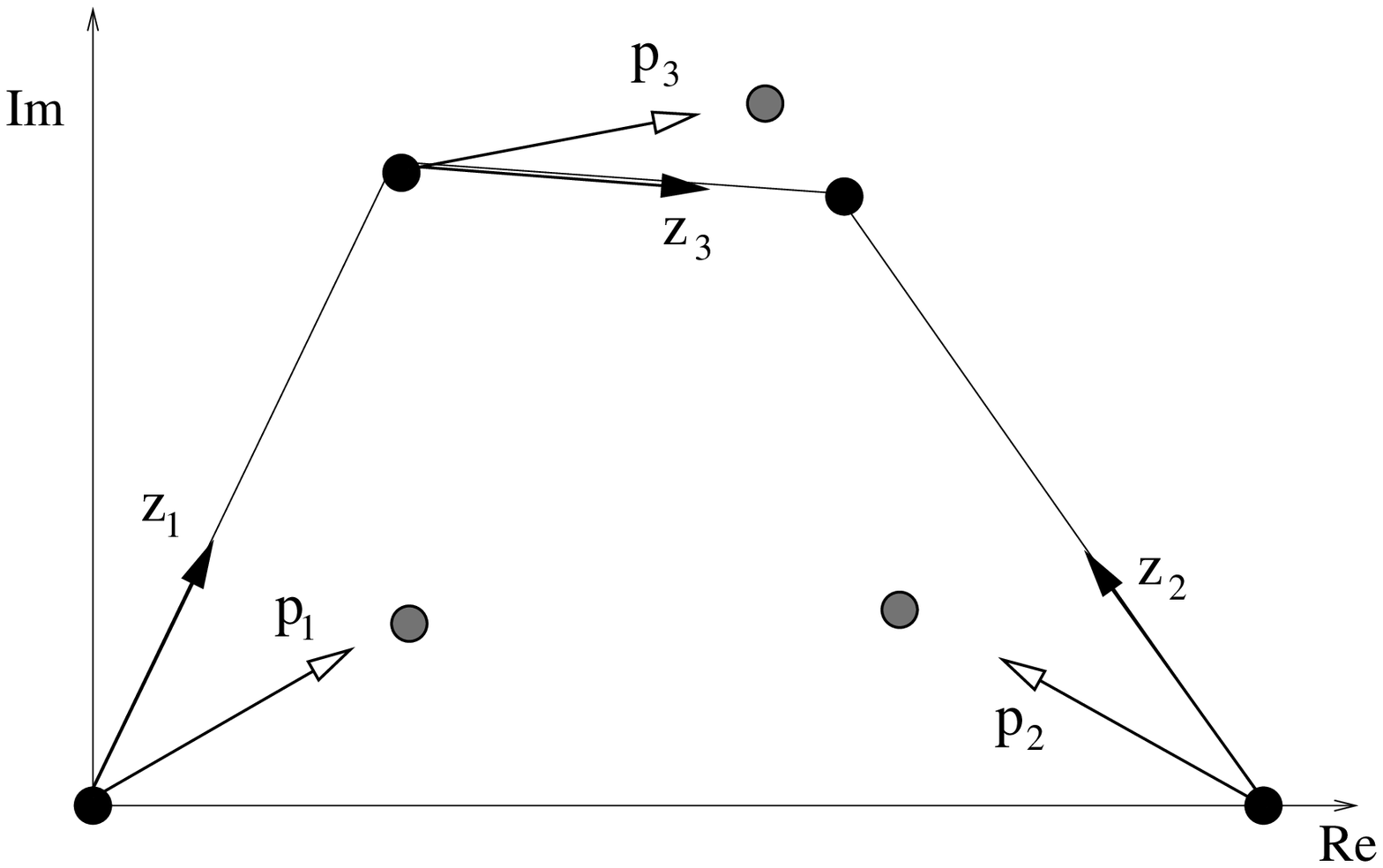} &
 	\includegraphics[width=5.5cm]{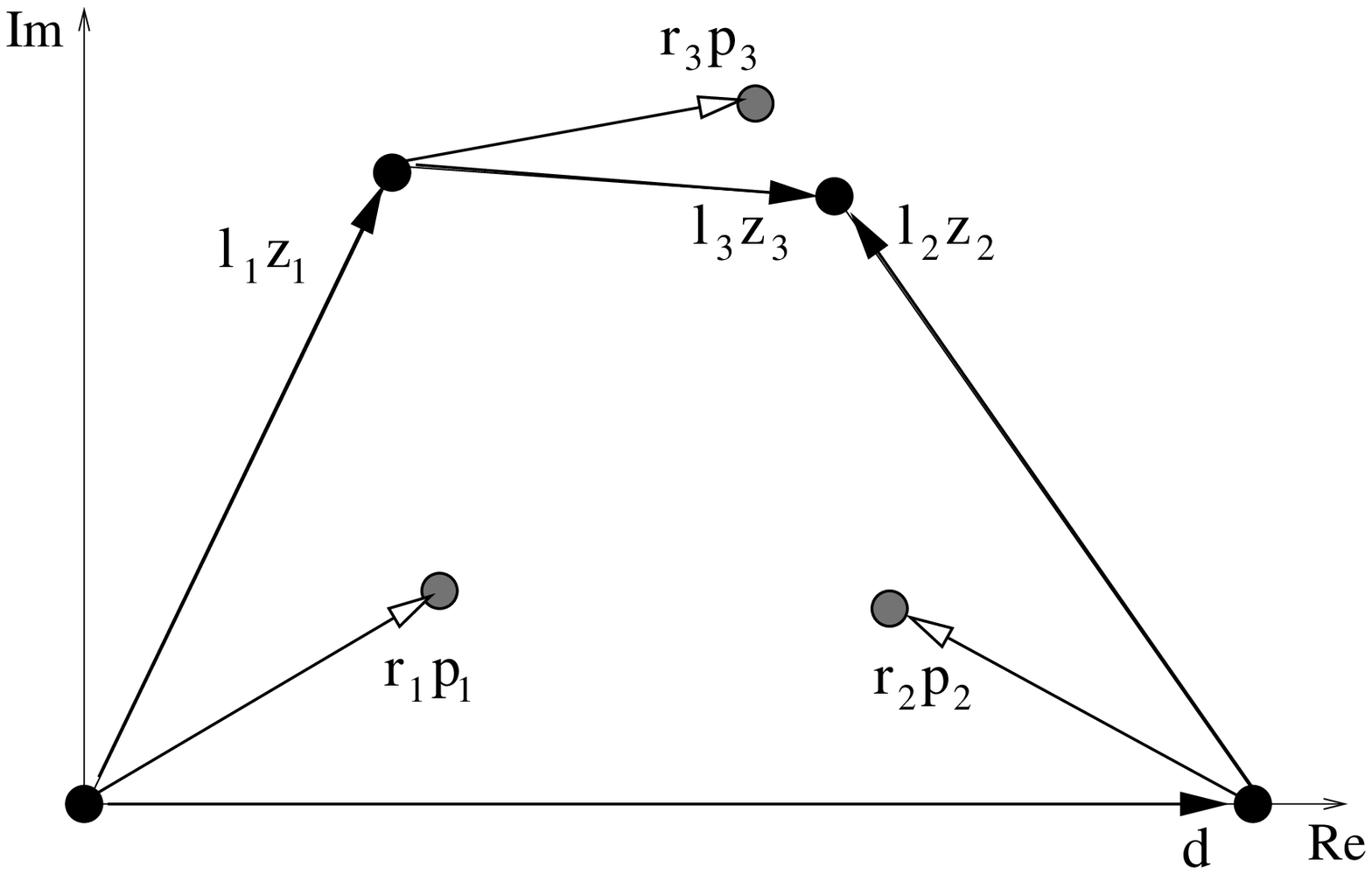} 
\end{tabular} 
\caption{Complex representation for the centre of masses.}
\label{fig:bar4Representation2}
\end{figure}

The links are connected by revolute joints rotating about axes pointing in a direction
orthogonal to the plane of motion.
The joint angles are specified using the time variables $\theta_1(t), \theta_2(t)$ and $\theta_3(t)$ as
shown in Figure (\ref{fig:bar4}).
Since the mechanism has only one degree of freedom, 
there is a relationship between these joint angles, which will be described below.
The kinematics of planar mechanisms can be conveniently represented in the complex plane, 
using complex numbers to describe the mechanism's configuration (Figure \ref{fig:bar4Representation})
and the location of the centre of mass (Figure \ref{fig:bar4Representation2}).
Referring to Figures \ref{fig:bar4Representation} and \ref{fig:bar4Representation2}, let
$\z_1,\z_2,\z_3$ be time dependent unit complex numbers and $\p_1, \p_2, \p_3$ unit complex numbers depending on the design parameters (actually only on $\psi_1, \psi_2, \psi_3$).
The orientation of $\p_i$ is specified relative to $\z_i$, i.e., it is attached to $\z_i$ and moves 
with it. If $\p_i$ coincides with $\z_i$, then $\p_i=1$. 


\subsection{Kinematic model}

The dependency between the different joint angles is described by the following closure constraint:
\begin{equation}
 \z_3 = \GName_1 \z_1 + \GName_2 \z_2 + \GName_3
\label{eq:LoopEqPos}
\end{equation}

\noindent
where $\GName_1, \GName_2, \GName_3 \in \R$ with $\GName_1 = \frac{-l_1}{l_3}, \GName_2 = \frac{l_2}{l_3}, \GName_3 = \frac{d}{l_3}$.
Taking the time derivative of equation (\ref{eq:LoopEqPos}), we get a relationship between the joint
angular velocities $\dt{1},\dt{2}$ and $\dt{3}$, namely:

\begin{equation}
 \z_3 \dt{3} = \GName_1 \z_1 \dt{1} + \GName_2 \z_2 \dt{2}
\label{eq:LoopEqVel}
\end{equation}
 
Since $\z_3$ is a unit complex number, $\z_3 \overline{\z_3} = \z_3 \z_3^{-1} = 1$
and therefore we obtain the following geometric constraint:
\begin{equation}
\GName = \left( \GName_1 \z_1 + \GName_2 \z_2 + \GName_3 \right) \left( \GName_1 \z_1^{-1} + \GName_2 \z_2^{-1}
+ \GName_3 \right) - 1 = 0.
\label{eq:GeoConstraint}
\end{equation}

The time derivative of the geometric constraint (\ref{eq:GeoConstraint}) can be written as a  
linear combination of the joint angular velocities:
\begin{equation}
 i (\KName_1 \dt{1} + \KName_2 \dt{2}) = 0
\label{eq:KinConstraint}
\end{equation}
where 
\begin{equation}
 \KName_1 = \GName_1 \GName_2 (\z_1 \z_2^{-1} - \z_1^{-1} \z_2) + \GName_1 \GName_3 (\z_1 - \z_1^{-1})
\end{equation}

\begin{equation}
 \KName_2 = \GName_1 \GName_2 (\z_1^{-1} \z_2 - \z_1 \z_2^{-1}) + \GName_2 \GName_3 (\z_2 - \z_2^{-1})
\end{equation}
It is noted that since $\KName_1$ and $\KName_2$ are purely imaginary, only one constraint
equation is obtained, over the real set.

\subsection{Position of the centre of mass}

Let $M$ be the total mass of the mechanism ($M=m_1+m_2+m_3$). 
The centre of mass ${\bf C}$ is:

\begin{equation}
{\bf C} = \frac{1}{M} \left( \FNameCoeff_1 \z_1 + \FNameCoeff_2 \z_2 + \FNameCoeff_3 \right) 
\label{eq:COM}
\end{equation}

where $\FNameCoeff_1, \FNameCoeff_2, \FNameCoeff_3 \in \C$:
\begin{eqnarray}
\label{eq:defF1}\FNameCoeff_1 & = & m_1 r_1 \p_1 + m_3 l_1 + \GName_1 m_3 r_3 \p_3\\
\label{eq:defF2}\FNameCoeff_2 & = & m_2 r_2 \p_2 + \GName_2 m_3 r_3 \p_3\\
\label{eq:defF3}\FNameCoeff_3 & = & m_2 d + \GName_3 m_3 r_3 \p_3.
\end{eqnarray}

 These equations were derived in \cite{MooreSchichoGosselin:2007}\footnote{In this paper, indices 2 and 3 for the links have been permuted.}.

\subsection{Angular momentum of the mechanism}

\begin{figure}
\begin{tabular}{c}
	\includegraphics[width=9cm]{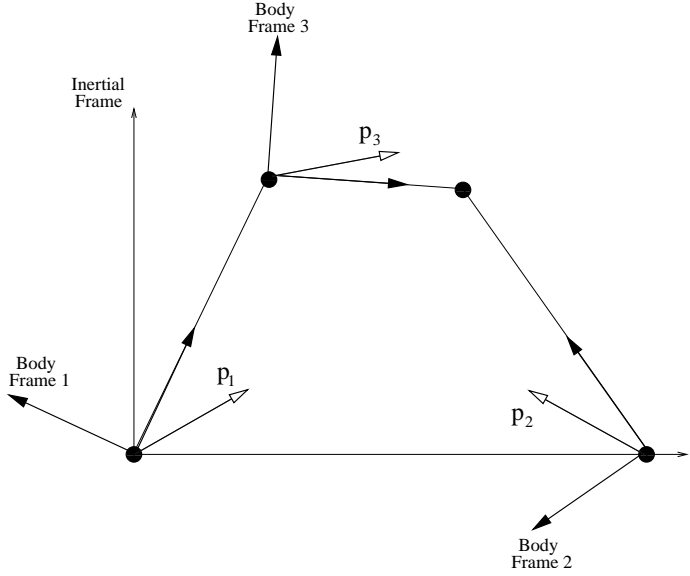}
\end{tabular} 
\caption{Unit vectors representation}
\label{fig:bar4Frames}
\end{figure}

Since the mechanism is planar, the contribution of body $i$ to the angular momentum is a scalar and can be given 
in the following form:
\begin{equation}
\HName_i = m_i \left\langle \overrightarrow{r}_{i/0} , -i \dot{\overrightarrow{r}}_{i/0} \right\rangle + I_i \dt{i}
\end{equation}
where $\overrightarrow{r}_{i/0}$ and $\dot{\overrightarrow{r}}_{i/0}$ are respectively the position 
and the velocity of the centre of mass of body $i$ with respect to a given inertial frame, 
$I_i$ denotes the moment of inertia of body $i$ with 
respect to its centre of mass and $\left\langle \ast, \ast \right\rangle$ is the scalar product of planar vectors.
The total angular momentum $H$ of the system is given by the sum of the angular momentum of the links, 
(i.e. $\HName = \HName_1 + \HName_2 + \HName_3$).

The angular momentum of the first body with respect to the inertial frame is:
\begin{equation}
 \HName_1 = \langle r_1 \p_1 \z_1, - i m_1 (i r_1 \p_1 \z_1 \dt{1}) \rangle + I_1 \dt{1} =  \\
\langle r_1 \p_1 \z_1, m_1 r_1 \p_1 \z_1 \dt{1} \rangle + I_1 \dt{1} = J_1 \dt{1} 
\label{eq:H1}
\end{equation}

The contribution of the second body to the angular momentum is given by:
\begin{equation}
 \HName_2 = \left\langle \left( d + r_2 \p_2 \z_2 \right), m_2 r_2 \p_2 \z_2 \dt{2} \right\rangle 
+ I_2 \dt{2} = \left[ m_2 d r_2 \left(\frac{\p_2 \z_2 + \p_2^{-1} \z_2^{-1}}{2}\right) + J_2\right] \dt{2}
\label{eq:H2}
\end{equation}

For the third body, we get 
\begin{equation}
 \HName_3 = \left\langle \left( l_1 \z_1 + r_3 \p_3 \z_3 \right), m_3 \left(l_1 \z_1 \dt{1} + r_3 \p_3 \z_3 \dt{3} \right) \right\rangle + I_3 \dt{3} \\
\label{eq:H3}
\end{equation}

Substituting equations (\ref{eq:LoopEqPos}) and (\ref{eq:LoopEqVel}) into
equation (\ref{eq:H3}), we can eliminate $\z_3 \dt{3}$ and $\dt{3}$ and obtain en expression in terms
of $\z_1, \z_2, \dt{1}, \dt{2}$ only.
The total angular momentum $H$ of the mechanism is then given by:
\begin{equation}
\HName = \HName_1 + \HName_2 + \HName_3 = \KName_3 \dt{1} + \KName_4 \dt{2}
\label{eq:AngMom}
\end{equation}
where $\KName_3$ and $\KName_4$ are written as:
\begin{eqnarray*}
\KName_3 & = & a_1 \z_1 + a_2 \z_1^{-1} + b_1 \z_1 \z_2^{-1} + b_2 \z_1^{-1} \z_2 + c\\
\KName_4 & = & u_1 \z_2 + u_2 \z_2^{-1} + v_1 \z_1 \z_2^{-1} + v_2 \z_1^{-1} \z_2 + w
\end{eqnarray*}
where constants $a_1, a_2, b_1, b_2, c, u_1, u_2, v_1, v_2$ and $w$ can be obtained from equations (\ref{eq:H1}, \ref{eq:H2}, \ref{eq:H3}, \ref{eq:AngMom}).

\subsection{Static and dynamic balancing}
In our settings, a mechanism is said to be statically balanced if the centre of mass of the mechanism remains 
stationary for infinitely many configurations (i.e. infinitely many choices of the joint angles).
From equation (\ref{eq:COM}), this condition can be formulated as:

\begin{equation}
 \FNameCoeff = \FNameCoeff_1 \z_1 + \FNameCoeff_2 \z_2 - C' = 0
\label{eq:StatBal}
\end{equation}
where $C'= C M - \FNameCoeff_3$ is a constant.
A mechanism is said to be dynamically balanced if the total angular momentum remains constant
for any motion of the mechanism.
In other words, the reaction forces and torques at its base induced by its motion are identically 
equal to 0, at all times.
Clearly, the mechanism must be statically balanced and the angular momentum must be constant.
Since the centre of mass is fixed and we want static balancing for all possible configurations and all joint
angular velocities, the angular momentum should therefore be 0, i.e.,
\begin{equation}
H=\KName_3 \dt{1} + \KName_4 \dt{2}=0 \label{eq:AngMom=0}
\end{equation}
Therefore, equations (\ref{eq:GeoConstraint}, \ref{eq:KinConstraint}, \ref{eq:StatBal}, \ref{eq:AngMom=0}) have to be satisfied. 
Among these four equations, only two (equation \ref{eq:KinConstraint} and \ref{eq:AngMom=0}) depend (linearly) on the joint
angular velocities and they can be rewritten in the following form:
\begin{equation}
\left[ \begin{array}{cc}
 \KName_1 & \KName_2 \\
 \KName_3 & \KName_4
\end{array} \right]
\left[ \begin{array}{c}
 \dt{1} \\
 \dt{2} 
\end{array} \right]
= \left[ \begin{array}{c}
 0 \\
 0 
\end{array} \right]
\label{eq:KinMatrixEq}
\end{equation} 
If the rank of the matrix $A =  \left[ \begin{array}{cc}
 \KName_1 & \KName_2 \\
 \KName_3 & \KName_4
\end{array} \right]$ is 2 and since the system is homogeneous, then the only solution is $\dt{1}=\dt{2}=0$. 
In other words, the mechanism is not moving. 
Therefore we must have:
\begin{equation}
 \KName := det(A)= \KName_1 \KName_4 - \KName_2 \KName_3 = 0
\label{eq:DynBal}
\end{equation}

Using this equation together with the geometric and static balancing constraints,
the joint angular velocities ($\dt{1}, \dt{2}$) are eliminated 
and a system of three algebraic equations (equation \ref{eq:GeoConstraint},\ref{eq:StatBal},\ref{eq:DynBal}) in terms of the unit complex variable $\z_1,\z_2$
is obtained.
In section 4, all possible design parameters satisfying this system of equations for infinitely many configurations
will be derived. However, we first need to introduce concepts and tools from toric geometry in order to solve
these equations. This is the subject of section 3.


\section{Toric geometry}

\begin{Definition}
A {\em Laurent polynomial} g over a ring $R$ is a formal sum
of monomials 
$x^\alpha := c_\alpha x_1^{\alpha_1} x_2^{\alpha_2}$,
where $x=(x_1,x_2)$ is a fixed pair of variables, and $\alpha\in\Z^2$,
$c_\alpha\in R$.
Its {\em support} is the set of all $\alpha \in \Z^2$ 
with non-zero coefficients $c_\alpha$.
Its {\em Newton polygon} is the convex hull of the support in $\R^2$.
The Laurent polynomials form a ring, namely
$R[x_1,x_2,(x_1x_2)^{-1}]$.
\end{Definition}

\begin{Definition}
The {\em Minkowski sum} of two convex sets $A$ and $B \subset \R^n$ is defined as:
\begin{equation}
A + B = \left\lbrace a + b \mid a \in A \wedge b \in B \right\rbrace 
\end{equation}
Note that $A+B$ is also a convex set.
\end{Definition}

\begin{Theorem}
\label{th:factor}
Assume that $R$ does not have zero divisors.
If $F, G$ are two Laurent polynomials, then:
\begin{equation}
\NewtonPolygon{FG} = \NewtonPolygon{F} + \NewtonPolygon{G}
\end{equation}
\end{Theorem}

We refer to Ostrowski, \cite{Ostrowski:1921, Ostrowski:1975} for a proof.

\begin{Remark} 
The assumption that $R$ has no zero divisor can be replaced by the weaker assumption that the corner coefficients of $G$, i.e., the coefficients at the vertices of $\NewtonPolygon{G}$, are not zero divisors.
\end{Remark}

In order to find out whether a given polynomial $G$ divides another given polynomial $F$, we introduce Laurent polynomial division.

\begin{Definition} 
Assume that $G$ is a Laurent polynomial such that its corner coefficients are no zero divisors. A finite subset $\Gamma$ of $\Z^2$ is called a {\em remainder support set} with respect to $G$ iff no multiple of $G$, except zero, has support contained in $\Gamma$.
\end{Definition}

\begin{Definition} 
\label{df:div}
Let $G$ be a Laurent polynomial such that its corner coefficients are invertible in $R$. Let $F$ be an arbitrary Laurent polynomial. Then $(Q,R)$ is a {\em quotient remainder pair} for $(F,G)$ iff the following conditions are fulfilled.

a) $F=QG+R$.

b) The support of $R$ is contained in $\NewtonPolygon{F}$.

c) The support of $R$ is a remainder support set with respect to $G$.
\end{Definition}

Quotient remainder pairs are not unique.
Here is a nondeterministic algorithm that computes quotient remainder pairs.

\begin{algorithm}
\caption{Toric Polynomial Division Algorithm}
\label{alg:ToricPolyDiv}
\begin{algorithmic}[1]
\STATE {\bf Input:} $F$, $G$, such that $G$ has invertible corner coefficients.
\WHILE{$F \ne 0$}
\STATE Select a linear functional $h:\R^2\to\R$, $(x,y)\mapsto(\alpha x+\beta y)$ such that $\alpha/\beta$ is irrational.
\STATE Compute the point $f\in\supp{F}$ that maximizes $h$.
\STATE Compute the point $g\in\supp{G}$ that maximizes $h$.
\IF {$\NewtonPolygon{G}+f-g\subset\NewtonPolygon{F}$} 
\STATE $M:=\frac{\coef{F,f}}{\coef{G,g}}x^{f_1-g_1}y^{f_2-g_2}$;
\STATE $Q:=Q+M$; $F:=F-MG$;
\ELSE
\STATE $M:=\coef{F,f}x^{f_1}y^{f_2}$;
\STATE $R:=R+M$; $F:=F-M$;
\ENDIF
\ENDWHILE
\STATE {\bf Output:} $Q$, $R$.
\end{algorithmic}
\end{algorithm}

\begin{Theorem}
Algorithm~\ref{alg:ToricPolyDiv} is correct.
\end{Theorem}

\begin{Proof}
The Newton polygon of $F$ becomes smaller in each {\bf while} loop, 
hence it is clear that Algorithm~\ref{alg:ToricPolyDiv} terminates. 
Also, any monomial which is added to $R$ is contained in $\NewtonPolygon{F}$,
hence it follows that $R$ fulfills (b) in Definition~\ref{df:div}.

No step in the algorithm changes the value of $F+QG+R$. Initially,
this value is the given polynomial $F$, and in the end, this value
is equal to $QG+R$. This shows that (a) in Definition~\ref{df:div}
is fulfilled.

In order to prove (c) in Definition~\ref{df:div}, we claim that the
following is true throughout the execution of the algorithm: if $H$ is
any Laurent polynomial such that $GH$ has support in 
$\NewtonPolygon{F}\cup\supp{R}$, then the coefficients of $GH$
at the exponent vectors in $\supp{R}$ are zero.

Initially, $\supp{R}$ is empty and the claim is trivially true. 
If the claim is true before step~8, then it is also true after step~8,
because this step does not change $R$ and does not increase the
Newton polygon of $F$. 

Assume that for a certain Laurent polynomial $H$, the claim is true 
before step~11 and false after step~11. Then it follows that the
coefficient of $GH$ at $f$ is not zero, because this is the only exponent
vector which is new in $R$. The support of $GH$ is also contained in
the Newton polygon of $F$ before step~11, hence $f$ is the unique
vector in $\supp{GH}$ where $h$ reaches maximal value. Because $g$
is the unique vector in $\supp{G}$ where $h$ reaches a maximal value,
it follows that $(f-g)\in\supp{H}$. Then 
$\NewtonPolygon{G}+f-g\subset\NewtonPolygon{GH}$ as a consequence
of Theorem~\ref{th:factor}. But this implies that the {\bf if} condition
in step~6 is fulfilled for $G$ and $F$ before step~11, and therefore
step~11 is not reached for such values of $F$ and $G$.

It follows that the claim is true throughout the execution of 
Algorithm~\ref{alg:ToricPolyDiv}. In particular, it is true at the end,
which shows that (c) in Definition~\ref{df:div} holds.
\end{Proof}

\begin{Example}

Let $F = c_{02} y^2 + c_{11} x y + c_{01} y + c_{10} x + c_{00}$ and
$G = d_{01} y + d_{10} x + d_{00}$. 
The result of the polynomial division algorithm is shown in table (\ref{table:DivisionExample}).
Therefore, $F$ is divisible by $G$ if and only if $R=0$, or in other words if all coefficients of $R$ are zero:
\begin{equation}
c_{11} - \frac{c_{02} d_{10}}{d_{01}}= c_{10} - \frac{A d_{10}}{d_{01}} =  c_{00} - \frac{A d_{00}}{d_{01}} = 0
\end{equation}

\begin{table}
\begin{tabular}{|c|c|c|c|}
\hline
Step & f & g & Computation \\
\hline \hline
7:  & (0,2) & (0,1) & $M = \frac{c_{02}}{d_{01}} y$   \\
8:  &       &       & $Q = \frac{c_{02}}{d_{01}} y$   \\
    &       &       & $F = \left( c_{11} - \frac{c_{02} d_{10}}{d_{01}}\right) x y +
                        \left( c_{01} - \frac{c_{02} d_{00}}{d_{01}}\right) y + c_{10} x + c_{00}$   \\
\hline
10: & (1,1) & (0,1) & $M = \left( c_{11} - \frac{c_{02} d_{10}}{d_{01}} \right) x y$  \\
11: &       &       & $R = \left( c_{11} - \frac{c_{02} d_{10}}{d_{01}} \right) x y$  \\
    &       &       & $F = \left( c_{01} - \frac{c_{02} d_{00}}{d_{01}}\right) y + c_{10} x + c_{00}$ \\
\hline
 7: & (0,1) & (0,1) & $M = \frac{1}{d_{01}} \left( c_{01} - \frac{c_{02} d_{00}}{d_{01}}\right) =: \frac{A}{d_{01}}$  \\
 8: &       &       & $Q = \frac{c_{02}}{d_{01}} y + \frac{A}{d_{01}} $  \\
    &       &       & $F = \left( c_{10} - \frac{A d_{10}}{d_{01}} \right) x + \left( c_{00} - \frac{A d_{00}}{d_{01}} \right)$  \\
\hline
10: & (1,0) & (0,1) & $M = \left( c_{10} - \frac{A d_{10}}{d_{01}} \right) x$  \\
11: &       &       & $R = \left(c_{11} - \frac{c_{02} d_{10}}{d_{01}}\right) x y + \left( c_{10} - \frac{A d_{10}}{d_{01}} \right) x$  \\
    &       &       & $F = c_{00} - \frac{A d_{00}}{d_{01}} $  \\
\hline
10: & (0,0) & (0,1) & $M = c_{00} - \frac{A d_{00}}{d_{01}}$  \\
11: &       &       & $R = \left(c_{11} - \frac{c_{02} d_{10}}{d_{01}}\right) x y + \left( c_{10} - \frac{A d_{10}}{d_{01}} \right) x + \left(c_{00} - \frac{A d_{00}}{d_{01}} \right)$  \\
    &       &       & $F = 0$ \\
\hline
\end{tabular}
\caption{Example: Toric polynomial division algorithm}
\label{table:DivisionExample}
\end{table}

\end{Example}

\section{Balancing}

\subsection{Problem description}
\begin{figure}
	\begin{tabular}{ccc}
		\includegraphics[height=3.5cm]{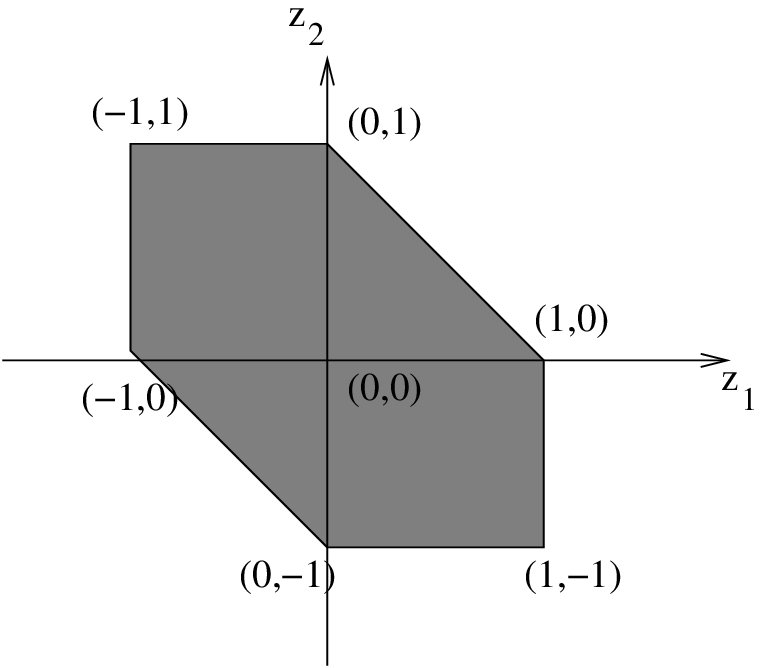}&
		\includegraphics[height=3.5cm]{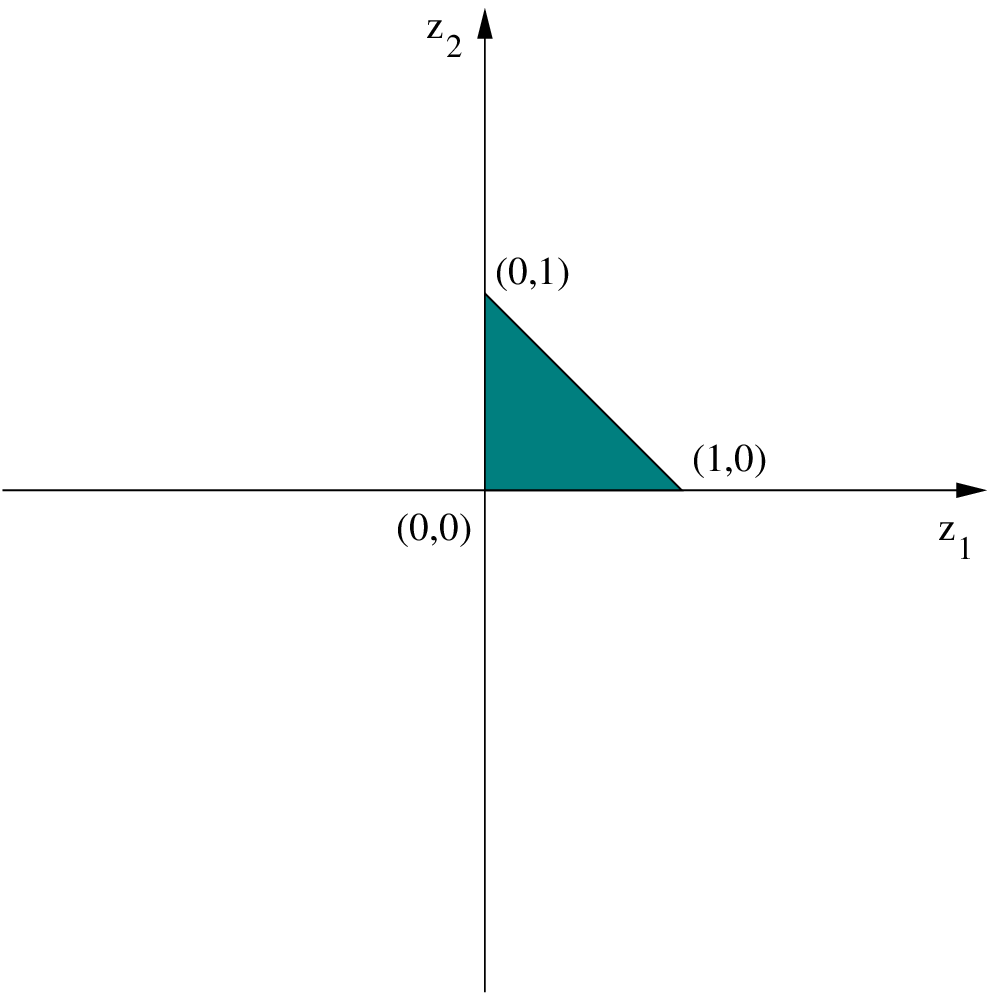} \hspace{1cm}&
		\includegraphics[height=3.5cm]{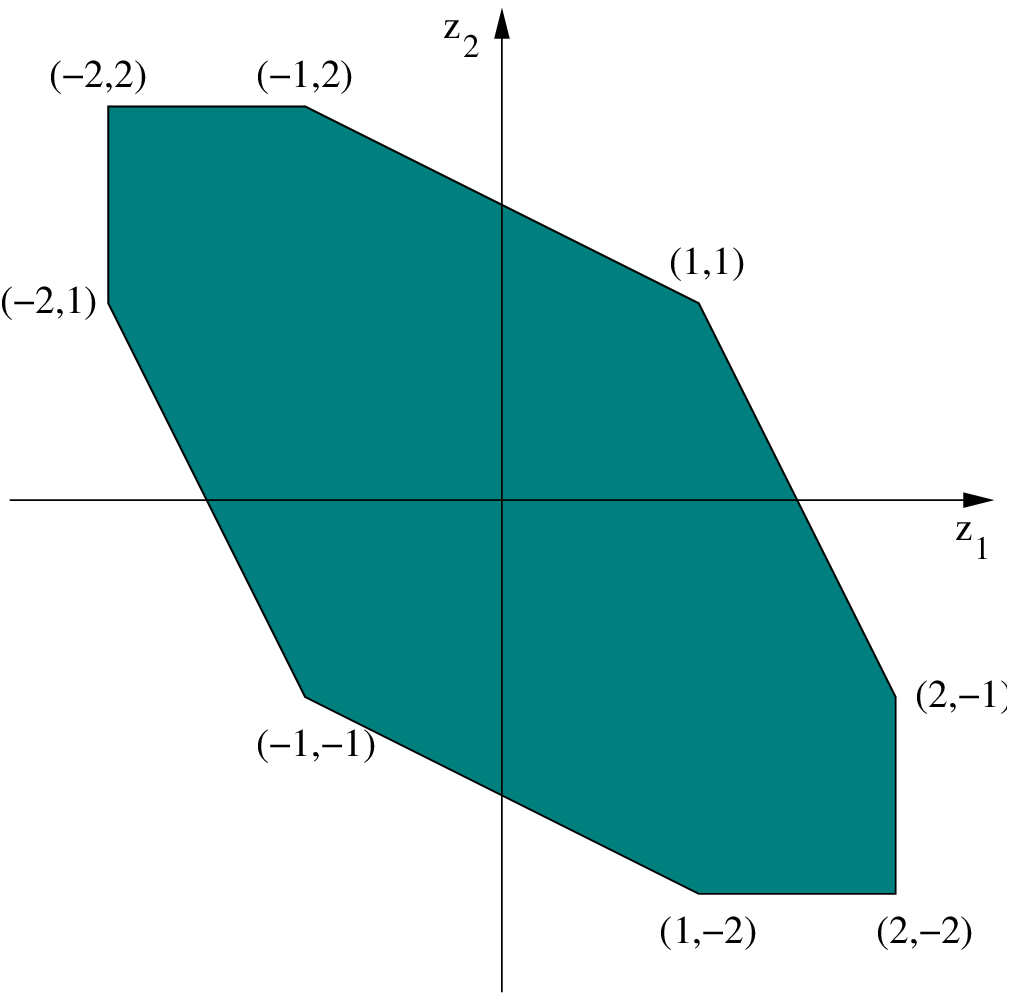}\\
                $\NP{\GName}$&$\NP{\FName}$&$\NP{\KName}$
	\end{tabular}
	\caption{Newton polygons of $\GName$, $\FName$ and $\KName$.}
	\label{fig:gfNewtonPolygon}
\end{figure}

The problem addressed in this paper can be stated as follows:
find all possible design parameters such that there exists a valid non-constant trajectory of the
planar four-bar mechanism for which the mechanism is dynamically balanced (i.e.: it is statically balanced 
and the angular momentum of the system is 0).
Formally, let $\KName = \left\lbrace (\z_1, \z_2 \in \C^2 \mid \GName(\z_1,\z_2) = 0) \right\rbrace$ be an
infinite set representing a valid non-constant trajectory.
For this trajectory, we want the mechanism to be statically balanced, i.e.:
\begin{equation}
\forall_{(\z_1,\z_2) \in \KName} \GName(\z_1,\z_2)=0  \Longrightarrow \FName(\z_1,\z_2)=0 
\label{eq:StaticBalancingFormulation}
\end{equation}
and the angular momentum to be 0:
\begin{equation}
\forall_{(\z_1,\z_2) \in K} \GName(\z_1,\z_2)=0   \Longrightarrow \KName(\z_1,\z_2)=0
\label{eq:DynamicBalancingFormulation}
\end{equation}
where $\GName$, $\FName$ and $\KName$ are defined in equations
(\ref{eq:GeoConstraint}, \ref{eq:StatBal}, \ref{eq:DynBal}).
Their Newton polygons are shown in figure \ref{fig:gfNewtonPolygon}.
Using the following theorem, we can reformulate this problem as a factorization problem of Laurent polynomials.

\begin{Theorem}
\label{th:IrrLaurentPoly}
Let $\GName$ be an irreducible Laurent polynomial.  Let $\KName$ be a Laurent polynomial(not necessarily irreducible).
Let $S \subseteq {\C^*}^2$ such that $\GName$ has infinitely many zeros in $S$.
The following are equivalent:
\begin{enumerate}
\item  $\forall (\z_1,\z_2) \in S, \GName(\z_1,\z_2)=0 \Rightarrow \KName(\z_1, \z_2) = 0$
\item $\exists$ Laurent polynomial $L$ such that $\KName = \GName \cdot L$
\end{enumerate}
\end{Theorem}

A proof of this theorem can be found in \cite{MooreSchichoGosselin:2007}.

\subsection{Static balancing}

Assume $\GName(\z_1,\z_2)$ is irreducible and using theorem (\ref{th:IrrLaurentPoly}) we are looking 
for a Laurent polynomials $L(\z_1,\z_2)$ such that

\begin{equation}
 \FName(\z_1,\z_2) = \GName(\z_1,\z_2) L(\z_1,\z_2)
\label{eq:factStat}
\end{equation}

If the geometric constraint $\GName$ is not irreducible, we consider all possible decomposition
of $\GName$ into irreducible components (see table \ref{table:gFact}). 
Every such decomposition imposes constraints on the kinematic parameters $l_1, l_2, l_3$ and $d$ (table \ref{table:gFact2}).
For a given decomposition, one component corresponds to a kinematic mode of the mechanical system.
For each of these decompositions and components, we can apply theorem (\ref{th:IrrLaurentPoly}).
Using this approach, necessary and sufficient conditions for the static balancing of planar
four-bar mechanisms can be obtained as shown in \cite{MooreSchichoGosselin:2007}.
These conditions are described in table \ref{table:Results}.

\begin{table}
\begin{tabular}{ccc}
&&\\
\includegraphics[height=1.3cm]{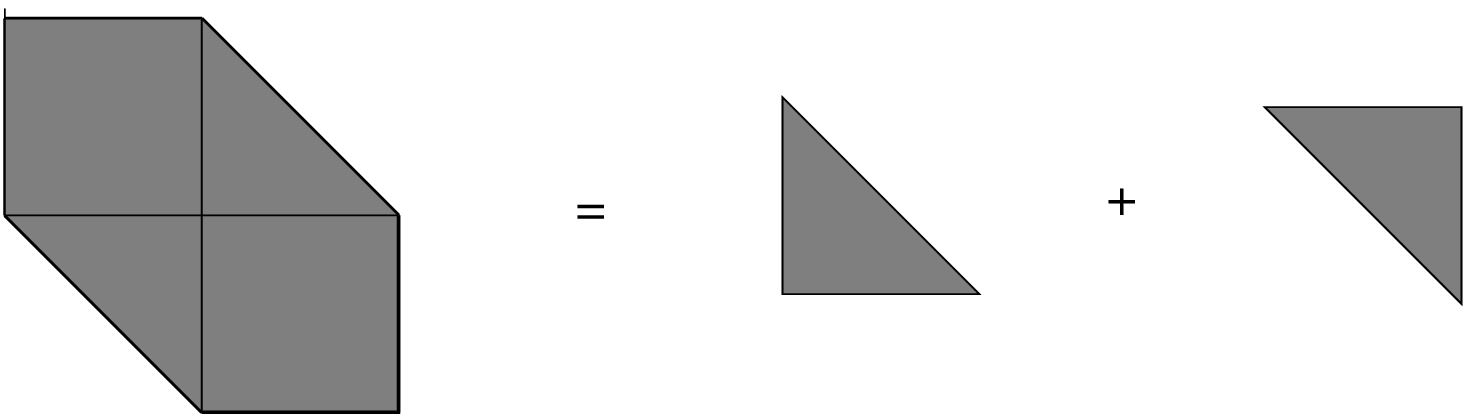} & \ \ \ \ \ \ &
\includegraphics[height=1.3cm]{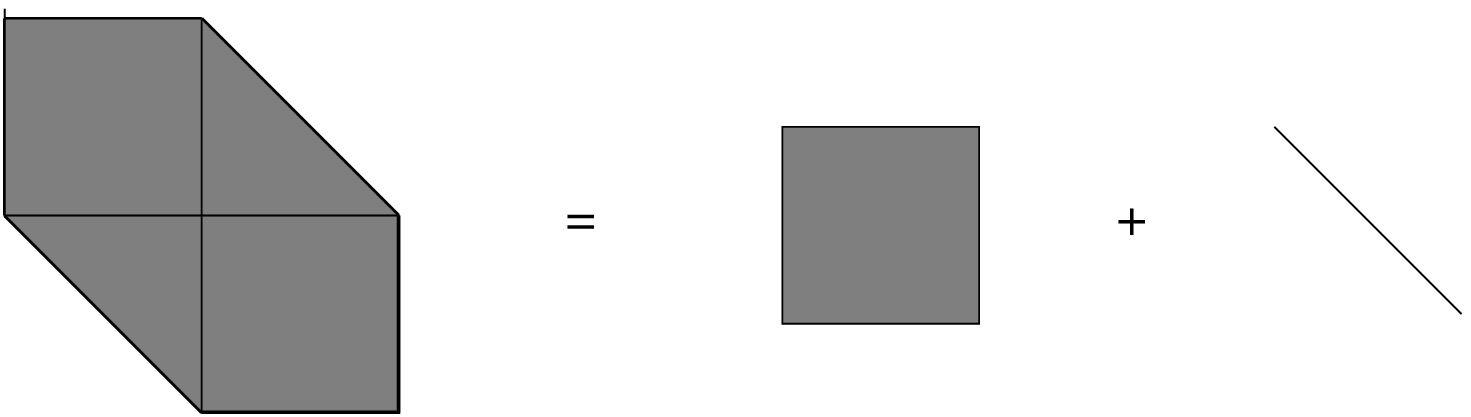}\\
(I)&&(II)\\
&&\\
&&\\
\includegraphics[height=1.3cm]{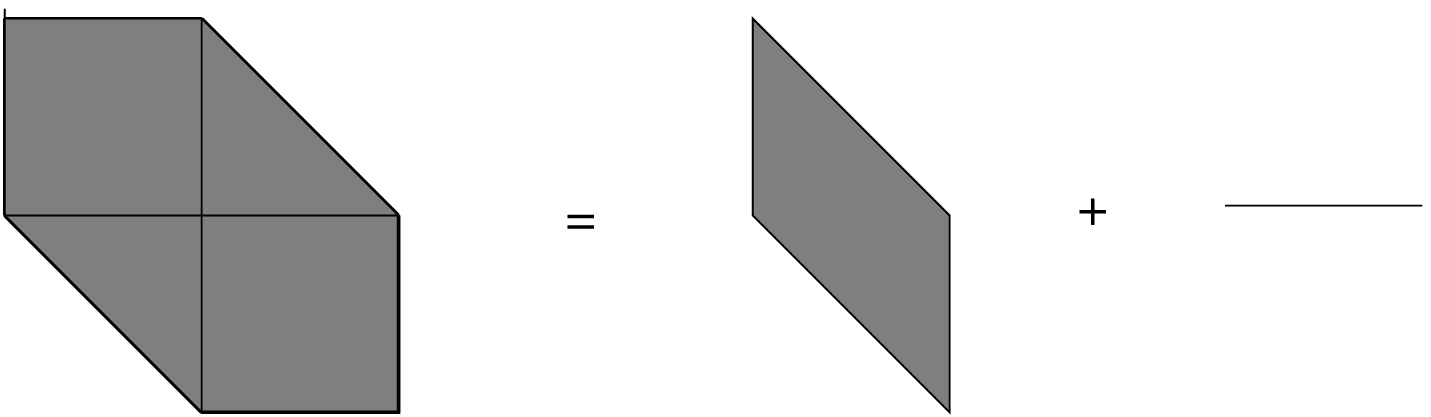}&&
\includegraphics[height=1.3cm]{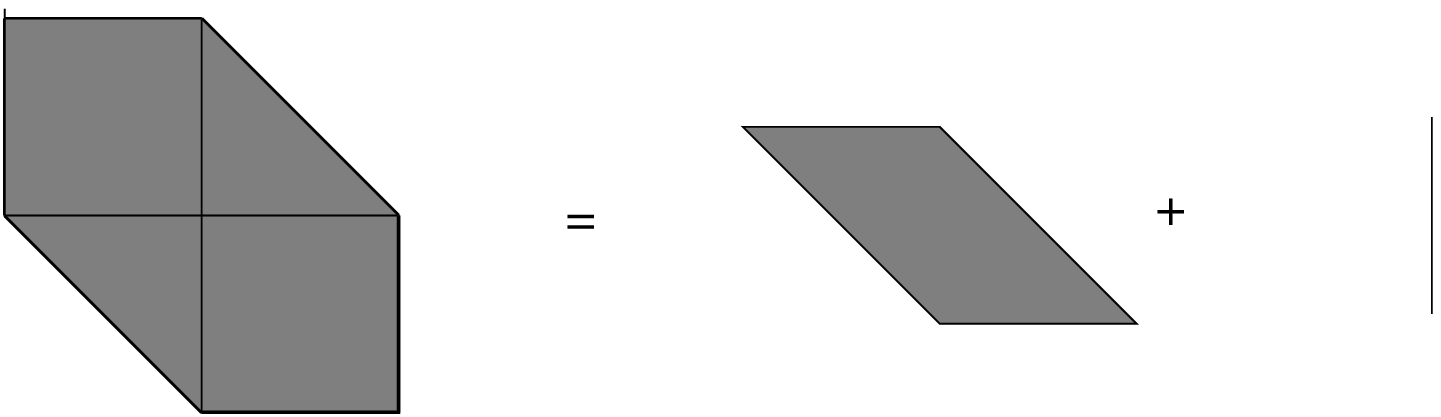}\\
(III)&&(IV)\\
&&\\
&&\\
\multicolumn{3}{c}{\includegraphics[height=1.3cm]{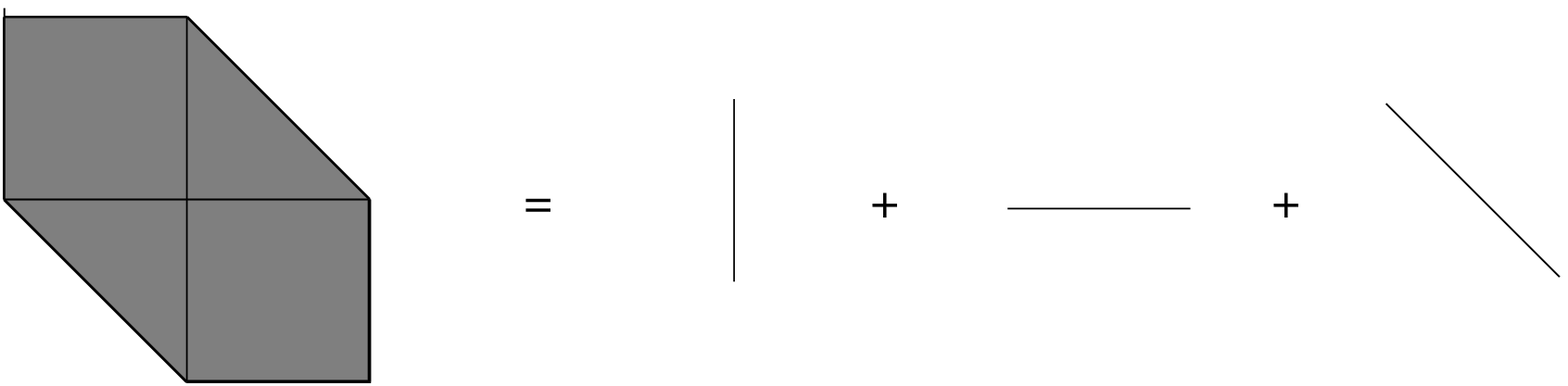}}\\
\multicolumn{3}{c}{(V)}\\
\end{tabular}
\caption{Possible decomposition of $\GName$ based on Newton polytopes and Minkowski sums.}
\label{table:gFact}
\end{table}

\begin{table}
\begin{tabular}{|c|c|c|c|c|}
\hline
Case & Kinematic constraint &   Mode A   & Mode B  &  Mode C  \\
\hline
&&&&\\ 
II&\includegraphics[width=30mm]{figure/gFact2MinkowskiSum.eps}   
&\includegraphics[width=20mm]{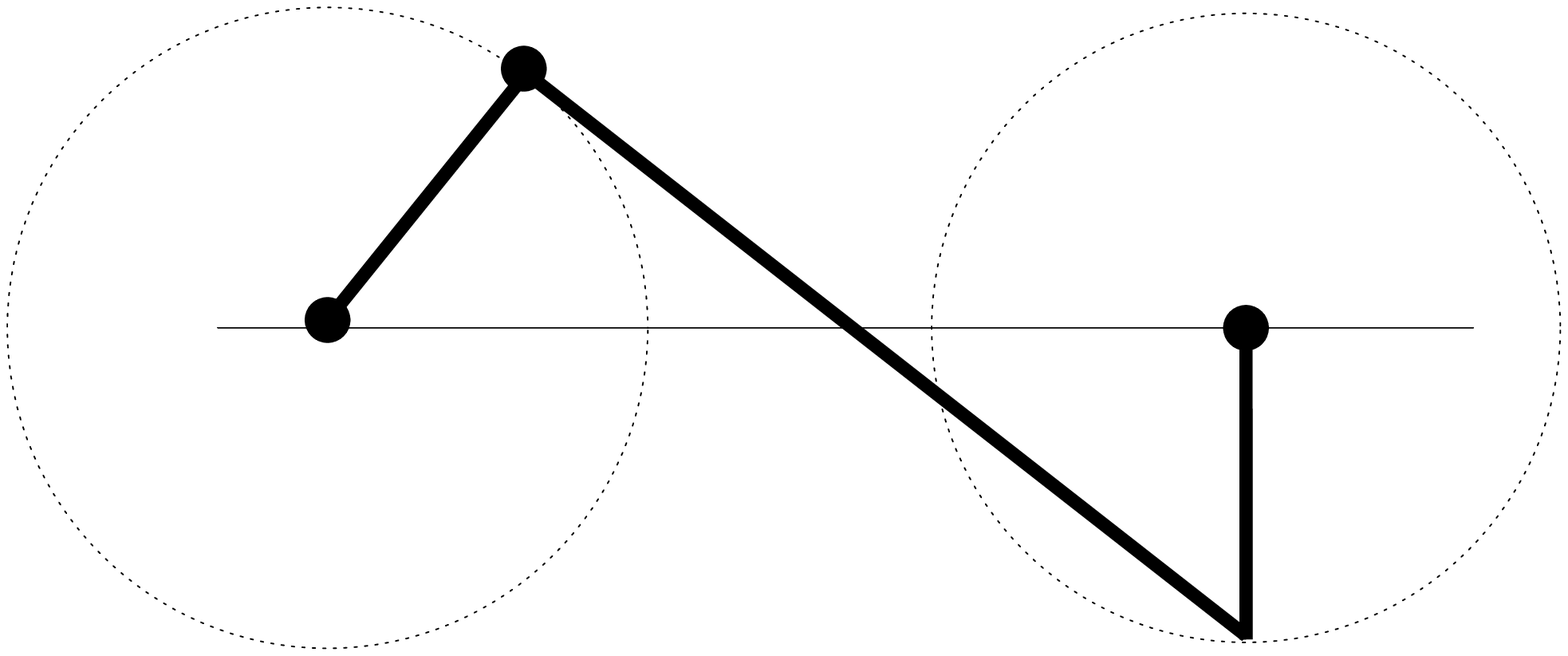}
& \includegraphics[width=20mm]{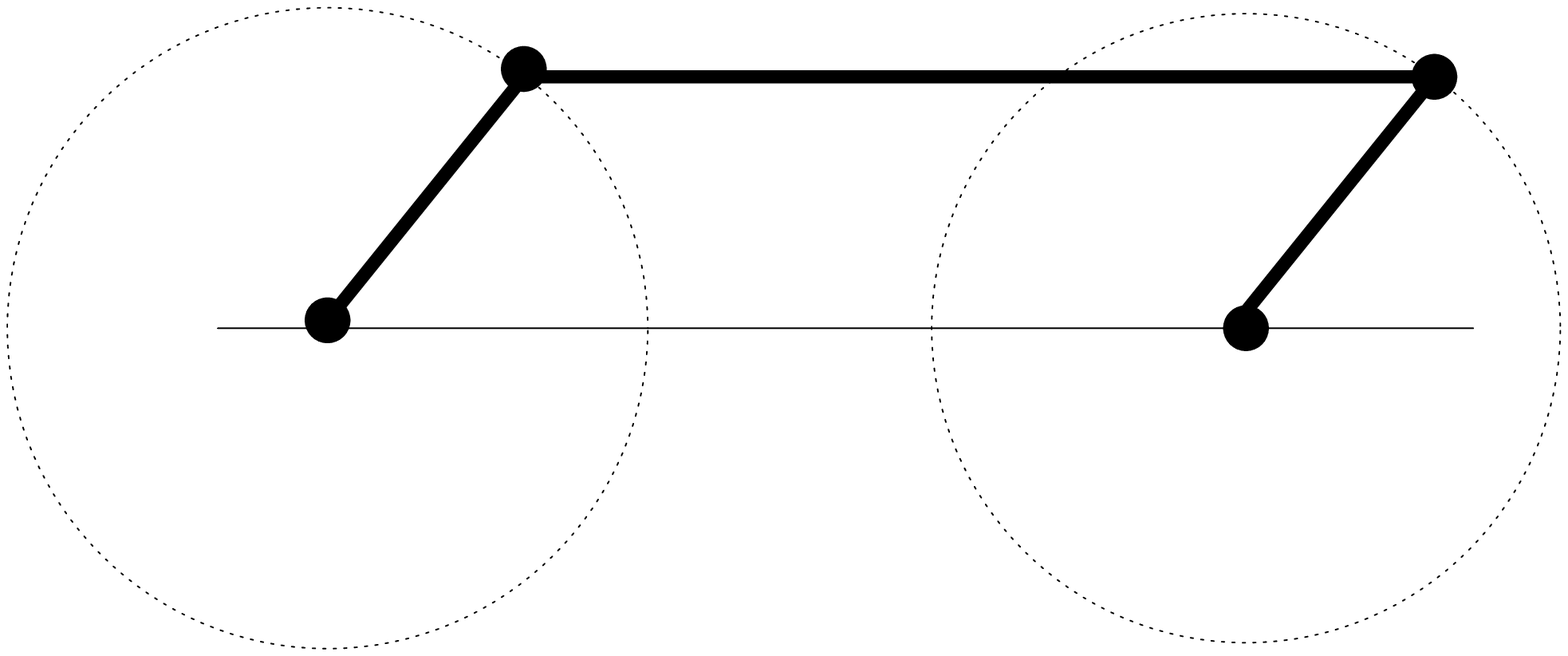} & \\
&$l_1 = l_2, l_3 = d, l_1 \neq l_3$&$\z_1 \neq \z_2$&$\z_1 = \z_2$&\\
&&&&\\
\hline
&&&&\\
III&\includegraphics[width=30mm]{figure/gFact4MinkowskiSum.eps}    
&\includegraphics[width=20mm]{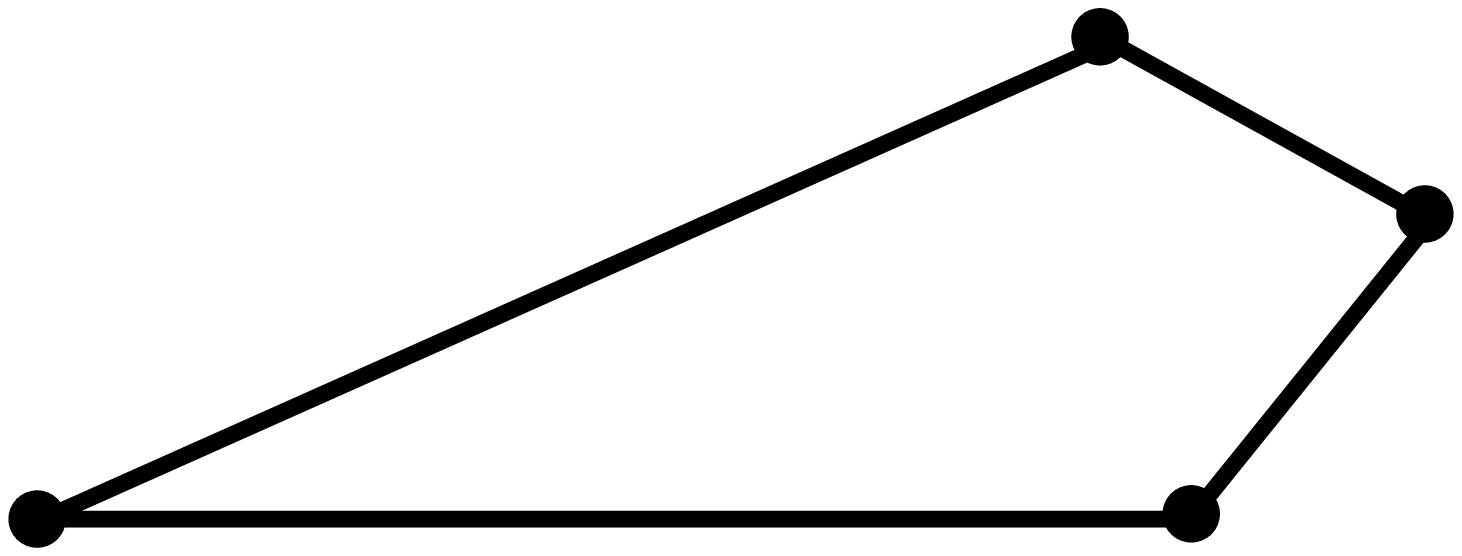}
&  \includegraphics[width=20mm]{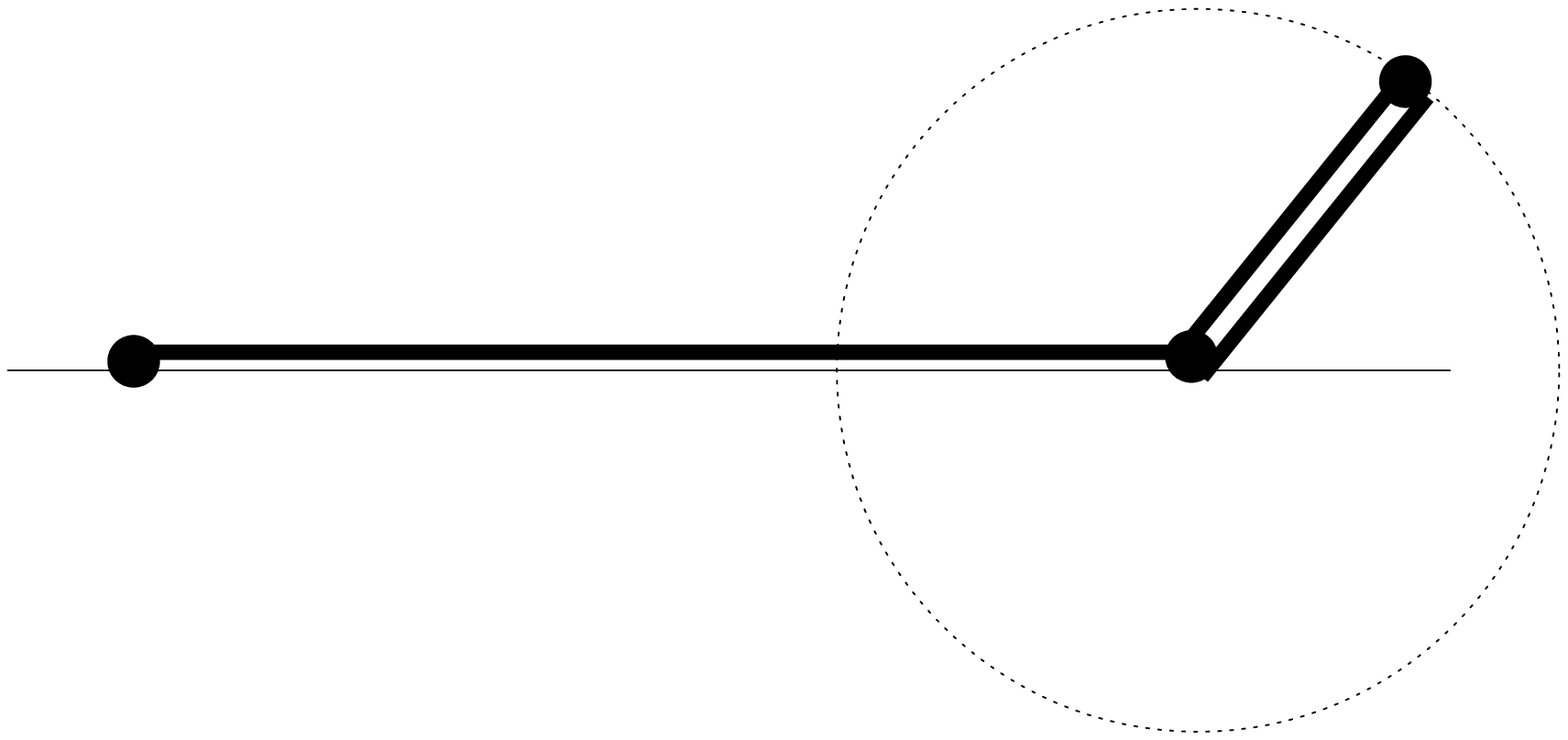} &\\
&$l_1 = d, l_2 = l_3, l_1 \neq l_2$&$\z_1 \neq 1$&$\z_1 = 1$&\\
&&&&\\
\hline
&&&&\\
IV&\includegraphics[width=30mm]{figure/gFact3MinkowskiSum.eps}  
&\includegraphics[width=20mm]{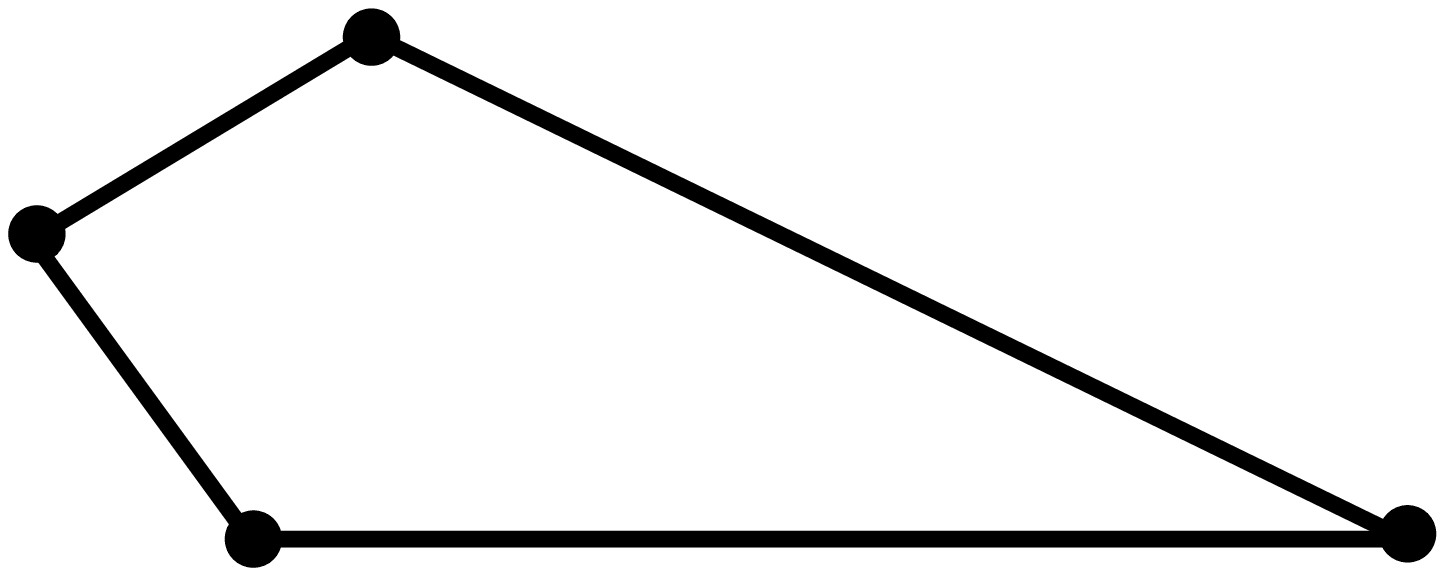}
& \includegraphics[width=20mm]{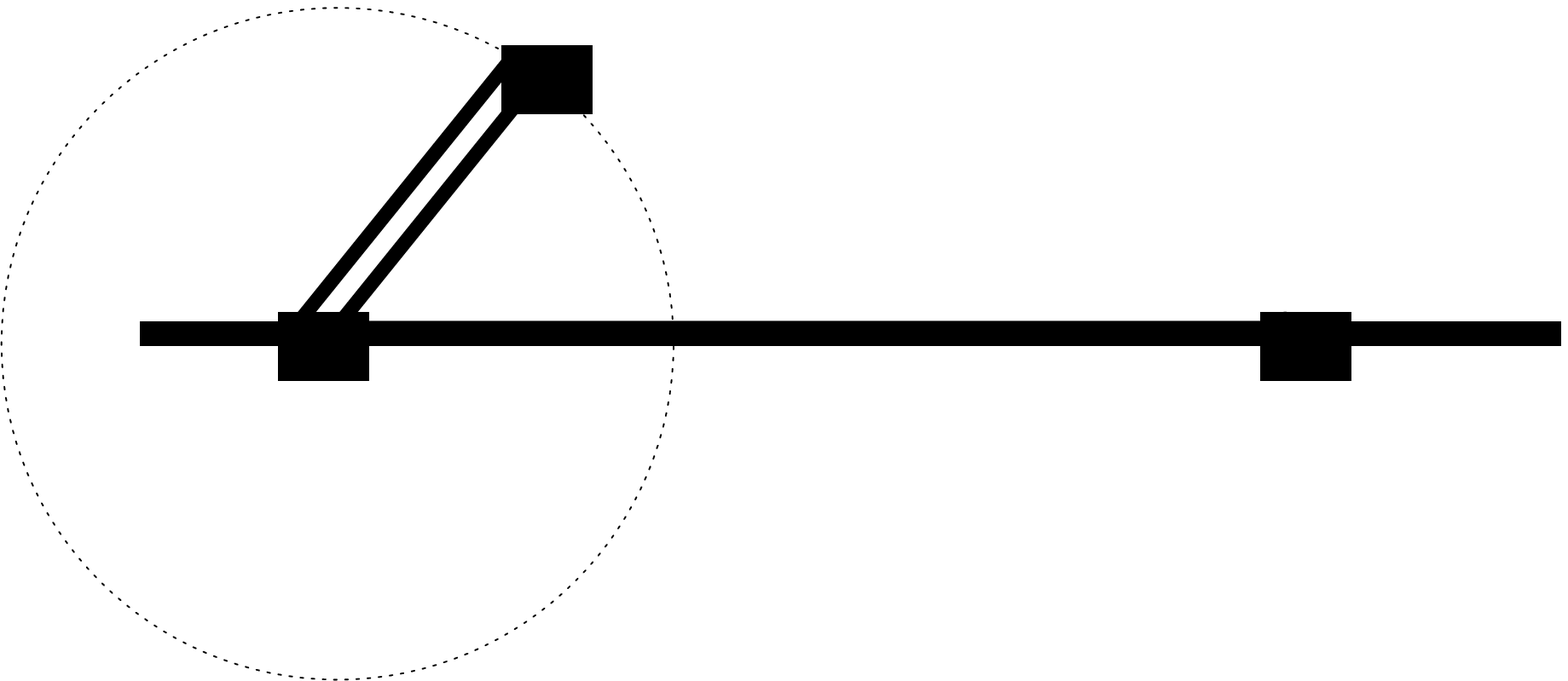} & \\
&$l_1 = l_3, l_2 = d, l_1 \neq l_2$&$\z_2 \neq 1$&$\z_2 = 1$&\\
&&&&\\
\hline
&&&&\\
V&\includegraphics[width=30mm]{figure/gFact5MinkowskiSum.eps}  
&\includegraphics[width=10mm]{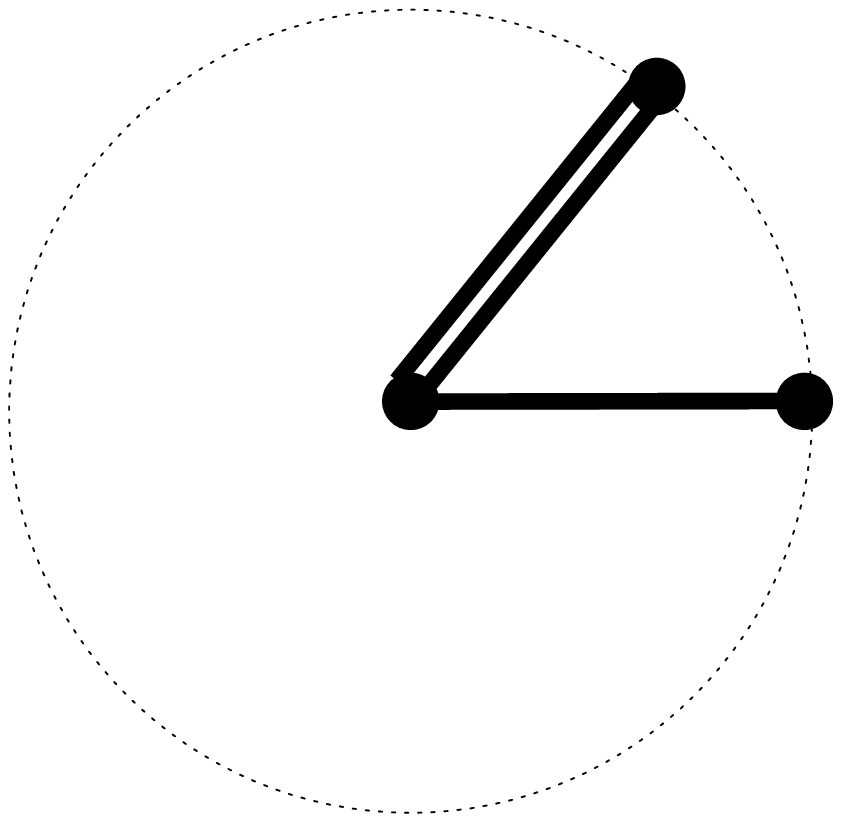} 
& \includegraphics[width=10mm]{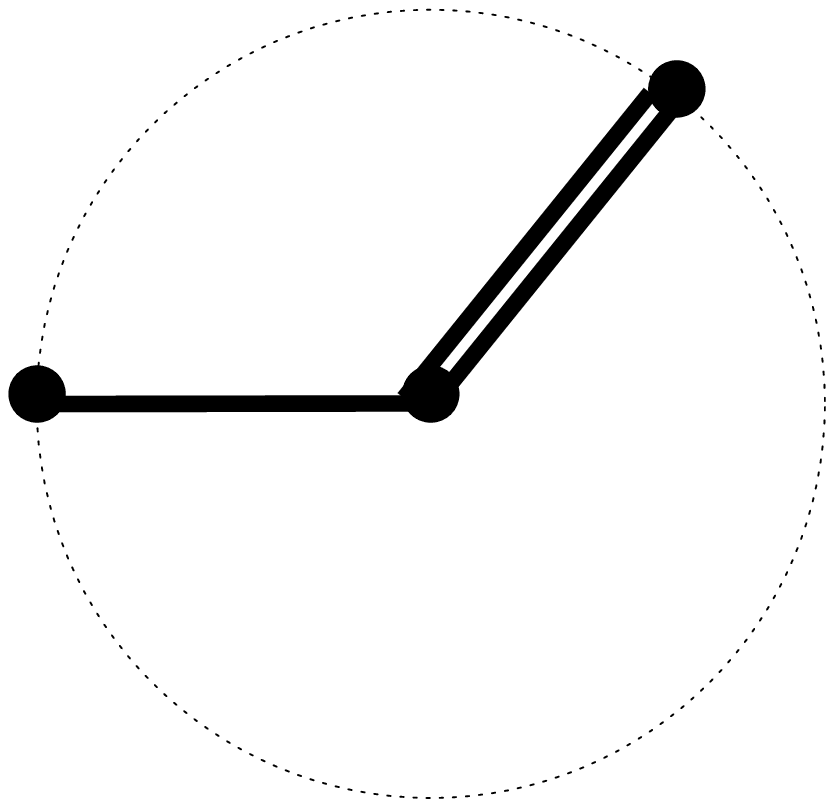} 
& \includegraphics[width=15mm]{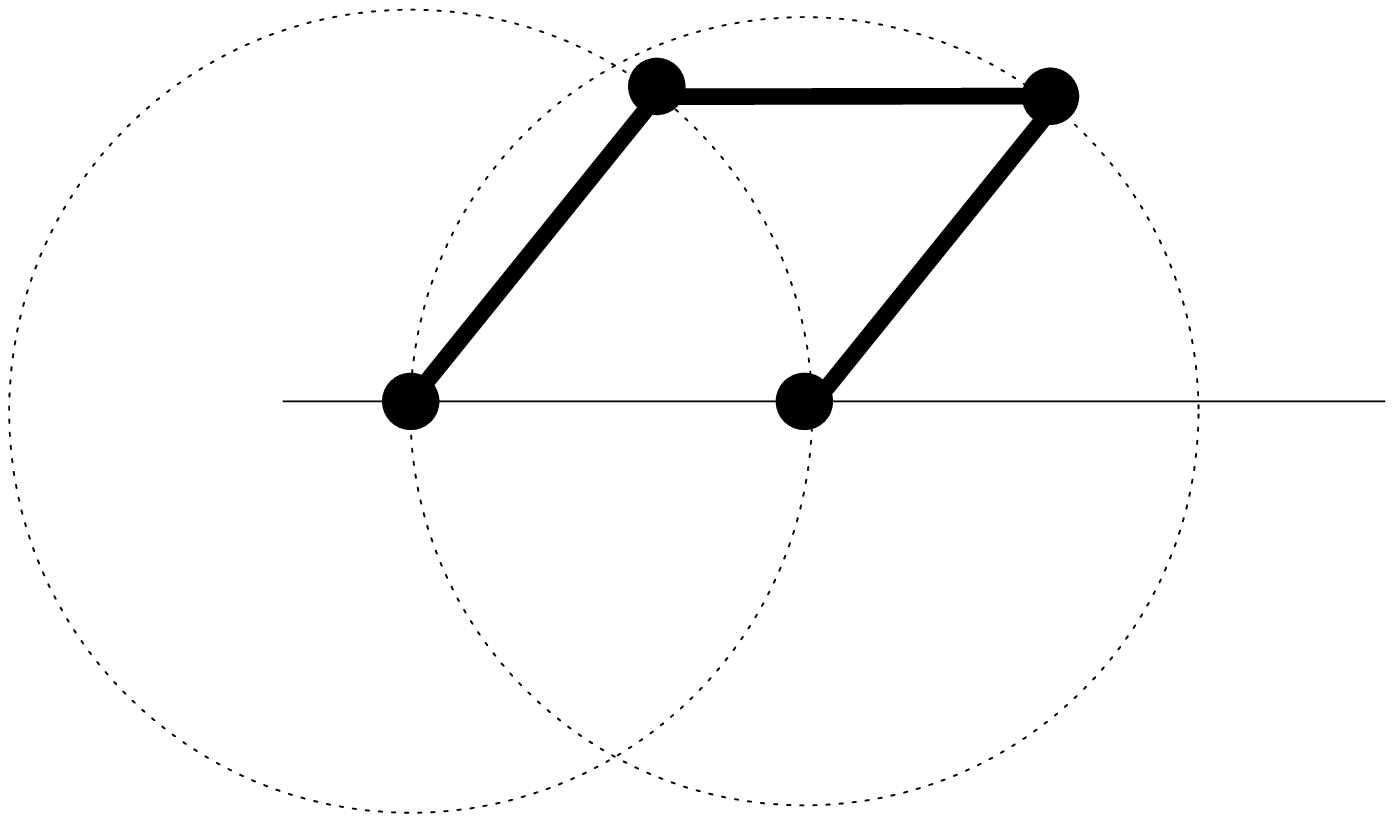}\\
&$l_1 = l_2 = l_3 = d$&$\z_2 = -1$&$\z_1 = 1$&$\z_1 = \z_2$\\
&&&&\\
\hline
\end{tabular}
\caption{Kinematic mode.}
\label{table:gFact2}
\end{table}

\subsection{Dynamic balancing}

The same approach as in the static balancing case could be used to find all sufficient and 
necessary condition for the dynamic balancing (equation \ref{eq:DynamicBalancingFormulation}).
However, due to the complexity of the Newton polygon of $\KName$, the approach would 
lead to a cumbersome case by case analysis, making it unpractical and prone to error. 
Using the toric polynomial division algorithm \ref{alg:ToricPolyDiv}, the same result can
be obtained in a semi-automatic way (using symbolic computation tools) without this case by case analysis.
In our case, the computations were performed using a Maple implementation of the Toric Polynomial Division algorithm and the solutions are presented below.

\subsubsection{Irreducible case}

Figure \ref{fig:PolDivCaseIrr} illustrates the main steps of the polynomial division algorithm 
when the geometric constraint $\GName$ is irreducible.
The algorithm gives a set of constraints in terms of the design parameters which can be combined with the
static balancing constraints.
Among these constraints, we obtain

\begin{equation}
 w = l_2^2 J_3 + l_3^2 J_2 = 0
\end{equation}

\noindent
with $J_2=m_2 r_2^2 + I_2$ and $J_3=m_3 r_3^2 + I_3$. Therefore $I_2=I_3=r_2=r_3=0$ which is physically not possible.
Therefore, if $\GName$ is irreducible, a planar four-bar mechanism cannot be dynamically balanced.

\begin{figure}
	\centering
        \psfrag{Act1}{PDiv 1}
        \psfrag{Act2}{PCst 1}
        \psfrag{Act3}{PCst 2}
        \psfrag{Act4}{PDiv 2}
        \psfrag{Act5}{PCst 3}
        \psfrag{Act6}{PCst 4}
        \psfrag{Act7}{PCst 5}
        \psfrag{Act8}{PCst 6}
        \psfrag{Act9}{PCst 7}
        \psfrag{ActA}{PDiv 3}
        \psfrag{ActB}{PCst 8}
        \psfrag{ActC}{PCst 9}
        \psfrag{ActD}{PCst 10}
        \psfrag{ActE}{PCst 11}
        \psfrag{ActF}{PCst 12}
        \psfrag{ActG}{End}
	\includegraphics[height=10cm]{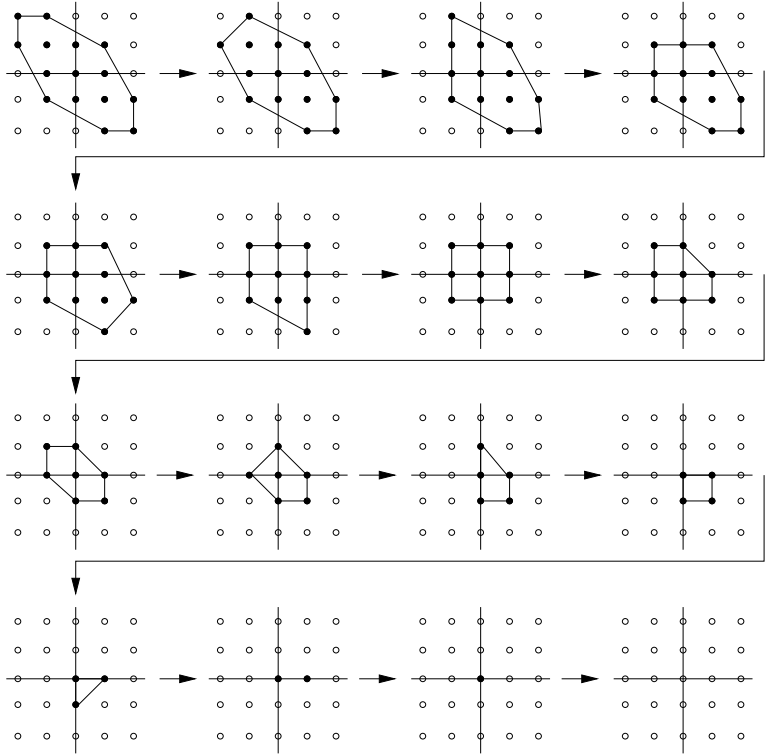}
	\caption{Polynomial division algorithm for the irreducible case}
	\label{fig:PolDivCaseIrr}
\end{figure}

\subsubsection{Reducible case II, Mode A}

Here we have $l_3=d$ and $l_1=l_2=:l$.
This is one case where we get solutions that are physically realizable.
In order to decribe the solutions, we introduce another set of
parameters, namely 
\[ q_i:= r_i p_i m_i , J_i:= I_i+m_ir_i^2 , i=1,2,3. \]
The variables $I_i,r_i,p_i$ can then be eliminated easily.
The balancing problems have an additive structure: if a balanced
mechanism picks up weights at each of its three bars in a balanced
way, then the composed mechanism is also dynamically balanced.
The variables are now chosen so that the balancing conditions
become linear in $m_i,q_i,J_i,i=1,2,3$. Here they are.
\begin{equation}
  q_1 = \frac{l}{d} q_3 - lm_3, q_2 = -\frac{l}{d} q_3 , 
\label{eq:cIIAq12}
\end{equation}
\begin{equation}
  J_1 = \frac{d^2+l^2}{d} q_3 - J_3 - l^2 m_3,
  J_2 = \frac{d^2-l^2}{d} q_3 - J_3 . 
\label{eq:cIIAJ12}
\end{equation}
A first consequence is that $q_1,q_2,q_3$ must be real.
The parameters fulfill also the inequality constraints
\begin{equation}
m_i>0, J_im_i - |q_i|^2 >0 
\label{eq:cIIAineqJ}
\end{equation}
for $i=1,2,3$. In particular, $J_1$ and $J_2$ must be positive.
By equation (\ref{eq:cIIAJ12}), we get an upper bound for $m_3$,
which must be larger than the lower bound from (\ref{eq:cIIAineqJ})
($i=3$). This yields
\begin{equation}
(dq_3-J_3)(dJ_3-l^2q_3) > 0 .
\label{eq:cIIAJ12a}
\end{equation}
It follows that $J_3$ is contained in the open interval 
$(dq_3,\frac{l^2}{d}q_3)$. 
(Note that $\frac{d^2-l^2}{d}q_3>0$ as a consequence of (\ref{eq:cIIAJ12}), 
that is why we know which of the two interval boundaries is bigger.) 
Then $q_3>0$ and $d>l$ follows. From $J_2>0$ and (\ref{eq:cIIAJ12}),
we get 
\begin{equation}
\frac{l^2}{d}q_3 < J_3 < \frac{d^2-l^2}{d}q_3 ,
\label{eq:cIIAJ12b}
\end{equation}
from which $d\ge \sqrt{2}l$ follows. 

Conversely, if $d\ge \sqrt{2}l$, then we can choose $q_3>0$ arbitrarily
and $J_3$ subject to (\ref{eq:cIIAJ12b}),
and $m_3$ between the upper and lower bound for $m_3$
derived above. Then (\ref{eq:cIIAJ12}) determines $J_1$ and $J_2$,
which will then be positive. Then (\ref{eq:cIIAq12}) determines $q_1$ and $q_2$,
and finally $m_1$ and $m_2$ can be chosen so that inequality (\ref{eq:cIIAineqJ})
is fulfilled. A possible solution is given in table (\ref{table:Example2A}).
\begin{table}
\begin{tabular}{|cccc|}
\hline
$l_1 = 1$              & $l_2 = 1$             & $l_3 = 4$           &    $d=4$   \\
$m_1 = m_1$            & $m_2 = \frac{m_1}{3}$ & $m_3 = \frac{2 m_1}{3}$           &            \\
$r_1 = \frac{1}{2}$    & $r_2 = \frac{1}{2}$   & $r_3 = 1$           &            \\
$p_1 = -1$               & $p_2 = -1$               & $p_3 = 1$           &            \\
$I_1 = \frac{3 m_1}{4}$  & $I_2 = \frac{5 m_1}{4}$  & $I_3 = \frac{m_1}{2}$ &            \\
\hline
\end{tabular}
\caption{Example of a dynamically balanced mechanism in case 2A}
\label{table:Example2A}
\end{table}


\subsubsection{Reducible case II, Mode B}
In this kinematic mode, the mechanism is a parallelogram with $z_1=z_2$.
Replacing $\z_1$ by $\z_2$, the constraint $\KName = 0$ becomes:
\begin{equation}
 -\frac{l_2}{d \z_2^2} (\z_2-1) (\z_2+1) \left[(a_1 + u_1) \z_2^2 + (b_1 + b_2 + c + v_1 + v_2 + w) \z_2 + u_2 + a_2 \right] = 0
\end{equation}

Therefore, the angular momentum vanishes if and only if one of these factors vanishes.
The first two factors (i.e.: $\z_2-1$ and $\z_2+1$) correspond to uncertainty
configurations in which it is possible to switch between mode A and mode B
(i.e.: where we can pass from one mode to the other) and are valid for only one configuration of the mechanism
($\z_1=1$ or $\z_1=-1$).
Therefore, the solution must come from the last factor and should be valid for all possible values of $\z_2$.
Therefore, the coefficients of $\z_2^2$, $\z_2$ and the constant term must vanish, i.e.:

\begin{equation}
 a_1 + u_1 = 0
\label{eq_solCase2B1}
\end{equation}
\begin{equation}
 b_1 + b_2 + c + v_1 + v_2 + w = 0
\label{eq_solCase2B2}
\end{equation}
\begin{equation}
 u_2 + a_2 = 0
\label{eq_solCase2B3}
\end{equation}
Equation (\ref{eq_solCase2B2}) can be written in terms of the design parameters in the following form:
\begin{equation}
 J_1 + l_2^2 m_3 + J_2 = 0
\end{equation}
This solution is physically not possible.

\subsubsection{Reducible case III, Mode A}
This case is completely symmetric with case IV mode A.

\subsubsection{Reducible case III, Mode B}

This case corresponds to the deltoid case with $z_1=1$.
By replacing $z_1=1$ in the equation of the angular momentum, $\KName$ can be written as:
\begin{equation}
\KName = \frac{d}{l_2 \z_2^2} (\z_2-1) (\z_2+1) ((u_1+v_2) \z_2^2 + w \z_2 + v_1 + u_2) = 0
\end{equation}

Using the same arguments as in case IIB, the term $(u_1+v_2) \z_2^2 + w \z_2 + v_1 + u_2 = 0$ 
should vanish for every unit complex number $z_2$, therefore:

\begin{equation}
 u_1+v_2 = 0
\label{eq_solCase3B1}
\end{equation}
\begin{equation}
 w = 0
\label{eq_solCase3B2}
\end{equation}
\begin{equation}
 v_1 + u_2 = 0
\label{eq_solCase3B3}
\end{equation}

Equation (\ref{eq_solCase3B2}) corresponds to $l_2^2 J_3 + l_3^2 J_2 = 0$ which is physically not possible.

\subsection{Reducible case IV, mode A}

Here we have $l_2=d$ and $l_1=l_3=:l$.
This is the second case where we get solutions that are physically realizable.
Again, we introduce the parameters $q_i$ and $J_i$ and eliminate
$I_i,r_i,p_i$ for $i=1,2,3$.
Here are the balancing conditions.
\begin{equation}
  q_1 = q_3 - lm_3, q_2 = -\frac{d}{l} q_3 , 
\label{eq:cIIIAq12}
\end{equation}
\begin{equation}
  J_1 = J_3 - l^2m_3,
  J_2 = \frac{l^2-d^2}{l} q_3 - J_3 . 
\label{eq:cIIIAJ12}
\end{equation}
It follows that $q_1,q_2,q_3$ must be real.
The parameters fulfill again the inequality constraints
\begin{equation}
m_i>0, J_im_i - |q_i|^2 >0 
\label{eq:cIIIAineqJ}
\end{equation}
for $i=1,2,3$; again it follows that $J_1$ and $J_2$ must be positive.
By equation (\ref{eq:cIIIAJ12}), we get an upper bound for $m_3$,
which must be larger than the lower bound from (\ref{eq:cIIIAineqJ})
($i=3$). This yields
\begin{equation}
(J_3-lq_3)(J_3+lq_3) > 0 .
\label{eq:cIIIAJn}
\end{equation}
This is equivalent to the statement $J_3>|lq_3|$. From (\ref{eq:cIIIAJ12}),
we obtain an upper bound for $J_3$, 
namely $\left|\frac{d^2-l^2}{l}q_3\right|$.
The lower bound must be larger than the upper bound, 
hence $q_3<0$ and $d\ge\sqrt{2}l$.

Conversely, if $d\ge \sqrt{2}l$, then we can choose $q_3<0$ arbitrarily and
$J_3$ between $-(lq_3)$ and $\frac{d^2-l^2}{l}q_3$.
Then we choose $m_3$ between the upper and lower bound for $m_3$
derived above. Next, (\ref{eq:cIIIAJ12}) determines $J_1$ and $J_2$,
which will then be positive. Then (\ref{eq:cIIIAq12}) 
determines $q_1$ and $q_2$,
and finally $m_1$ and $m_2$ can be chosen so that inequality 
(\ref{eq:cIIIAineqJ}) is fulfilled.

\subsubsection{Reducible case IV, Mode B}
This case is completely symmetric with case II1, mode B. 
Therefore, there are no solutions in this case.

\subsection{Reducible case V}
These cases are similar to II-B, III-B and IV-B. 
Using the same approach as above, the mechanism cannot be dynamically balanced.

\section{Conclusion}

The complete charaterization of dynamically balanced planar four-bar mechanisms is given in table \ref{table:Results}.
Note that these simple balanced mechanisms can be combined to build more complex 
planar and spatial balanced mechanism as shown in \cite{GosselinVollmerCoteWu:2004}.

\begin{table}
\begin{tabular}{|c|c|c|c|c|}
\hline
Kinematic constraints & Kinematic mode & Static balancing & Dynamic balancing \\
\hline
\hline
     $l_1=l_2, l_3=d$ & $\z_1 \neq \z_2$ & $\FNameCoeff_1=\FNameCoeff_2=0$  & possible iff $d \ge \sqrt{2} \  l_2$\\
                      & $\z_1=\z_2$      & $\FNameCoeff_1+\FNameCoeff_2=0$  & no \\
\hline
    $l_1=d, l_2=l_3$  & $\z_1 \neq 1$    & $\FNameCoeff_1=\FNameCoeff_2=0$  & possible iff $d \ge \sqrt{2} \  l_3$ \\
                      & $\z_1 = 1$       & $\FNameCoeff_2=0$           & no \\
\hline
    $l_2=d, l_1=l_3$  & $\z_2 \neq -1$   & $\FNameCoeff_1=\FNameCoeff_2=0$  & possible iff $d \ge \sqrt{2} \  l_3$ \\
                      & $\z_2 = -1$      & $\FNameCoeff_1=0$           & no \\
\hline
     Otherwise        &                  & $\FNameCoeff_1=\FNameCoeff_2=0$  & no \\ 
\hline
\end{tabular}
\caption{Complete balancing constraints for planar four-bar mechanisms.}
\label{table:Results}
\end{table}

\section*{Acknowledgement} 
The authors would like to thank Boris Mayer St-Onge
for the numerical validation of the examples using the software Adams.


\end{document}